\documentclass[review,hidelinks,onefignum,onetabnum]{siamart250211}

\usepackage{lipsum}
\usepackage{amsfonts}
\usepackage{amssymb}
\usepackage{graphicx}
\usepackage{epstopdf}
\usepackage{algorithmic}
\ifpdf
  \DeclareGraphicsExtensions{.eps,.pdf,.png,.jpg}
\else
  \DeclareGraphicsExtensions{.eps}
\fi
\usepackage{pgfplots}
\pgfplotsset{compat=1.18}
\usepgfplotslibrary{groupplots}
\usepackage{tikz}
\usepackage{multirow}
\usepackage{booktabs}
\usepackage{pgfplotstable}
\usepackage[justification=centering]{subcaption}
\usepackage{caption}
\usepackage[sort]{cite}
\usepackage{bm}
\newcommand{\tensor}[1]{\mathcal{#1}} 

\newsiamremark{remark}{Remark}
\newsiamremark{hypothesis}{Hypothesis}
\crefname{hypothesis}{Hypothesis}{Hypotheses}
\newsiamthm{claim}{Claim}
\newsiamremark{fact}{Fact}
\crefname{fact}{Fact}{Facts}

\headers{An Efficient Randomized Hierarchical Tucker Decomposition}{B. Leys, K. Meerbergen, N. Vannieuwenhoven, and M. Vivet}

\title{An Efficient Randomized Hierarchical Tucker Decomposition\thanks{
\funding{This work was in part funded by the Flanders Innovation \& Entrepreneurship (VLAIO) Baekeland mandate HBC.2025.0254 which is a collaboration between Siemens Industry Software and KU Leuven. N. V. was partially supported by the KU Leuven Interne Fondsen with project number C16/21/002, and by the Research Foundation -- Flanders (FWO) with projects G080822N and G011624N. M. V. was supported by the Flanders Innovation \& Entrepeneurship (VLAIO) O\&O Project Satisfy.AI (HBC.2025.0492).}}}

\author{Bram Leys\thanks{KU Leuven, Department of Computer Science, Research unit NUMA, B-3000 Leuven, Belgium. 
  ((\email{bram.leys@kuleuven.be},0009-0009-2844-3835); (\email{karl.meerbergen@kuleuven.be}, 0000-0002-1508-0248); (\email{nick.vannieuwenhoven@kuleuven.be}, 0000-0001-5692-4163).}
\and Karl Meerbergen\footnotemark[2]
\and Nick Vannieuwenhoven\footnotemark[2]
\and Mathijs Vivet\thanks{Siemens Digital Industries Software, Leuven, Belgium. 
	(\email{mathijs.vivet@siemens.com}, 0000-0003-2824-4034).}}

\usepackage{amsopn}

\ifpdf
\hypersetup{
 pdftitle={An Efficient Randomized Hierarchical Tucker Decomposition},
 pdfauthor={B. Leys, K. Meerbergen, N. Vannieuwenhoven, and M. Vivet}
}
\fi

\begin{document}
\usetikzlibrary{shapes.geometric, arrows.meta, positioning,calc,trees}
\maketitle

\begin{abstract}
We introduce randomized algorithms for computing the hierarchical Tucker decomposition of high-order tensors. Unlike traditional algorithms performed in a sequential leaves-to-root direction, our algorithms traverse the dimension tree in a root-to-leaves direction, leveraging randomized sequential projections to accelerate the computation of intermediate subspaces. This strategy can reduce the computation time by an order of magnitude relative to prior randomized algorithms using the leaves-to-root direction, while guaranteeing quasi-optimal approximation. Although the theoretical error bound is slightly larger than for leaves-to-root methods, in all our experiments the practical error adheres to the classic bound. Our results highlight the method's potential for scalable tensor decompositions in high-dimensional data analysis.
\end{abstract}


\section{Introduction}
Tensor decompositions are powerful tools for handling high-dimensional data organized in a multidimensional array, or tensor. Since many fields, ranging from computational physics to machine learning, generate and use such data~\cite{an_tensor_2019, bachmayr_low-rank_2023, richter_solving_2021}, tensor decompositions have diverse applications. They represent high-order multidimensional arrays using smaller, low-order tensors. Due to their reduced storage requirements, tensor decompositions also provide a framework for reduced-order or surrogate models of multivariate functions~\cite{sasagawa_construction_2023, kastian_two-stage_2020, strossner_approximation_2024, azaiez_data-driven_2021}.


This paper presents a more efficient algorithm for approximating a dense tensor by a \emph{hierarchical Tucker decomposition} (HTD) or tree tensor network~\cite{grasedyck_hierarchical_2010,hackbusch_new_2009}. The HTD represents the original tensor via a hierarchy of successively nested tensor product subspaces, in a way that refines the traditional Tucker decomposition~\cite{grasedyck_hierarchical_2010}. The HTD is parametrized by a tree graph, which defines the nesting of the subspaces. The tensor train decomposition (TT)~\cite{oseledets_tensor-train_2011} is an HTD where a line tree is used. 

Current algorithms to transform general dense tensors into the HTD format can suffer from high computational costs. These state-of-the-art algorithms perform a series of orthogonal projections onto subspaces defined by the dimension tree, traversing the tree and applying projections from the leaves to the root (LTR) to ensure the required nestedness property of HTDs. This paper, however, advocates for a \emph{root-to-leaves} (RTL) direction strategy. Our proposed algorithms are mathematically equivalent to an RTL method previously introduced by Hackbusch~\cite{hackbusch_tensor_2019}, who applied it to truncate tensors already represented as HTDs. While seemingly abandoned in favor of LTR traversal on theoretical complexity grounds as in, e.g.,~\cite{grasedyck_hierarchical_2010,hackbusch_tensor_2019, kressner_algorithm_2014,perisa_tensortoolboxjl_2021}, this article demonstrates that RTL traversal can accelerate HTD computation when \emph{randomized truncated singular value decompositions} are used, especially when sequential truncation is performed both across and within each layer. Experimental results demonstrate speedup factors of at least $8\times$ compared to the standard LTR algorithm of Grasedyck~\cite[Algorithm 2]{grasedyck_hierarchical_2010} for decomposing the Miranda computational fluid dynamics dataset~\cite{cook_tera-scalable_2005}. A detailed analysis of the approximation error and computational complexity is provided, along with experimental results that demonstrate the advantages of our approach. 

Our improvements build on insights from randomized numerical linear algebra, more specifically the calculation of singular value decompositions using randomized methods. Zniyed et al.~\cite{zniyed_tt-based_2020} already investigated the effectiveness such methods for the TT format, where it was shown that an exact decomposition could be sped up by following a more balanced tree-traversal strategy rather than the standard linear traversal of the TT-SVD algorithm~\cite{oseledets_tensor-train_2011}. This work extends these ideas to the HTD and broadens their applicability to general tensors, beyond those of exact rank by providing a theoretical quasi-optimal truncation error bound.
Related works applying randomized techniques to the HTD include~\cite{bucci_randomized_2024}, which employs randomized Nyström projections to accelerate tensor decomposition.
A comparison with this algorithm is made in \Cref{sec:experiments}, showing a speedup factor of one order of magnitude for our proposed methods.
Randomized sampling techniques also exist~\cite{ballani_black_2013,ryzhakov_black-box_2024} which view the data-generating process as a black box. This paper focuses on the dense tensor case.

The remainder of the paper is structured as follows: \Cref{sec:prelim} introduces the notation for tensors and their common operations used in this paper. \Cref{sec:htd} recalls the HTD and reviews related works. \Cref{sec:rtl} describes the new RTL algorithms and provides theoretical guarantees on the accuracy obtained when using them to truncate a general tensor into an HTD. These algorithms are applied to two realistic engineering datasets to evaluate their efficacy and generality in \Cref{sec:experiments}. Finally, \cref{sec:conclusions} gives some concluding remarks and discusses potential directions for further work.

\section{Preliminaries}
\label{sec:prelim}
We first fix notation and conventions for standard concepts and operations used throughout the paper. Unless stated otherwise, $\| . \|$, denotes the Frobenius norm. For a tensor $\tensor{A} \in \mathbb{R}^{n_1 \times \dots \times n_d}$, $\|\tensor{A}\| = \sqrt{\sum_{i_1 = 1}^{n_1} \dots \sum_{i_d = 1}^{n_d} |a_{i_1,\dots,i_d}|^2}$. The span of the columns of a matrix $A$ is denoted by $\operatorname{span}(A)$.

\subsection{Randomized singular value decomposition}
\label{sec:svd}
Let $A\in\mathbb{R}^{m\times n}$ be a matrix of rank $k$. 
The \emph{compact singular value decomposition (SVD)} factors $A$ as
\begin{equation*}
	A 
	= 
	U\,\Sigma\,V^\top
	=
	\begin{bmatrix}\widehat{U} & \widetilde{U}\end{bmatrix}
	\begin{bmatrix}\widehat{\Sigma} & \\ & \widetilde{\Sigma}\end{bmatrix}
	\begin{bmatrix}\widehat{V}^\top \\ \widetilde{V}^\top\end{bmatrix},
\end{equation*}
where $U\in\mathbb{R}^{m\times k}$ and $V\in\mathbb{R}^{n\times k}$ are matrices with orthonormal columns, $\Sigma \in\mathbb{R}^{k\times k}$ is a diagonal matrix with singular values $\sigma_1\geq\dots\geq\sigma_{k} > 0$ as elements, and $\widehat{U}$ and $\widehat{V}$ contain, respectively, the first $r$ left and right singular vectors corresponding to the $r$ largest singular values in the submatrix $\widehat{\Sigma}=\operatorname{diag}(\sigma_1,\dots,\sigma_r)$.

The \emph{rank-$r$ truncated SVD projection} is then $
\widehat{A} := \widehat{U}\,\widehat{\Sigma}\,\widehat{V}^\top$. By the Schmidt--Eckart--Young--Mirsky theorem~\cite{eckart_approximation_1936,mirsky_symmetric_1960, schmidt_zur_1907} this is the best rank-$r$ approximation of $A$ in (among others) the Frobenius norm:
\[
\|A-\widehat{A}\| = \min_{\operatorname{rank}(B)\le r}\|A-B\|.
\]
We call $\widehat{U}\widehat{U}^\top$ the \emph{rank-$r$ truncated SVD projector} since it satisfies $\widehat{U}\widehat{U}^\top A = \widehat{A}$.

Classic deterministic algorithms for computing the full SVD have a complexity of
$\mathcal{O}(m n \min(m,n))$~\cite{demmel_applied_1997}. Computing a truncated SVD by first computing a full SVD and then discarding $\widetilde{U}\widetilde{\Sigma}\widetilde{V}^\top$ becomes prohibitive for large-scale problems. 
\emph{Randomized SVD} methods, by contrast, approximate the leading singular vectors by projecting $A$ onto a low-dimensional random subspace~\cite{halko_finding_2011}. These methods reduce the cost to compute a rank-$r$ truncated SVD to $\mathcal{O}(m n r)$, giving large speedups if $r \ll \min(m,n)$. 

A particularly simple randomized SVD algorithm is a single sketching method with either a column space range finder or row space range finder based on the largest dimension of $A$~\cite{halko_finding_2011}. Assuming $m \ge n$, these sketching methods sample an $\Omega \in \mathbb{R}^{n\times r}$ from the Gaussian ensemble, i.e., all elements $\Omega_{i,j}$ are i.i.d. standard normal variables. This matrix is used to obtain an orthogonal basis of an approximation of the left dominant subspace of $A$ by \emph{sketching} $A$, i.e., computing $A\Omega$ and computing the QR decomposition $ A\Omega = QR$. Afterwards, a compact SVD is computed on the $r\times n$ matrix $Q^\top A = W\Sigma V^\top$ so that the rank-$r$ truncated SVD of $A \approx (QW) \Sigma V^\top$.
Despite their stochastic nature, their expected error can be minimized when the target rank is close to the rank of the matrix~\cite{halko_finding_2011}, enabling their use when an exact solution is unnecessary. Typically, a small oversampling constant, e.g., $p = 5$, is used to improve the accuracy. Then, the sketching uses $\Omega \in \mathbb{R}^{n\times (r + p)}$, while the compact SVD of $Q^\top A$ is truncated to rank $r$. More advanced techniques adaptively enlarge the orthogonal basis $Q$ until the expected approximation error falls under some threshold; see~\cite{halko_finding_2011}.

\subsection{Tensors}
\label{sec:tensors}
Throughout this article, $\mathrm{V}_1, \dots, \mathrm{V}_d$ denote finite-dimensional vector spaces over $\mathbb{R}$. 
A \emph{tensor $\tensor{A}$ of order $d$} is an element of the tensor product space \( \mathrm{V}_1 \otimes \cdots \otimes \mathrm{V}_d\), where $\otimes$ is the tensor product~\cite{greub_multilinear_1978}. If $n_k = \operatorname{dim}(\mathrm{V}_k)$, then the abstract space $\mathrm{V}_1 \otimes \dots \otimes \mathrm{V}_d$ can be identified with $\mathbb{R}^{n_1 \times \dots \times n_d}$, the space of real arrays~\cite{johnston_advanced_2021}.

A nonzero tensor $\tensor{B} \in \mathrm{V}_1 \otimes \cdots \otimes \mathrm{V}_d$ is called a \emph{rank-1 tensor} if it can be written as
\(\tensor{B} = \bm{v}_1 \otimes \cdots \otimes \bm{v}_d\) where each \(\bm{v}_k \in \mathrm{V}_k.\)
In coordinates such a tensor satisfies
\[
\tensor{B}_{i_1\ldots i_d}
= 
\bm{v}_{1,i_1}\,\bm{v}_{2,i_2}\cdots \bm{v}_{d,i_d},
\]
where $\bm{v}_{k,i_k}$ is the $i_k$-th entry of the vector $\bm{v}_{k}$.
Every tensor $ \tensor{A} \in \mathrm{V}_1 \otimes \cdots \otimes \mathrm{V}_d$ can, by definition~\cite{greub_multilinear_1978}, be expressed as a finite linear combination of rank-1 tensors:
\[
\tensor{A} = \sum_{i=1}^r \bm{v}_1^{(i)} \otimes \cdots \otimes \bm{v}_d^{(i)}, \quad \bm{v}_k^{(i)} \in \mathrm{V}_k.
\]

Let $\mathcal{I}_1 \sqcup \dots \sqcup \mathcal{I}_m = (1,\dots,d) =: \mathcal{I}$ be a partition of $m$ disjoint, nonempty, ordered tuples. 
We define the \emph{flattening} (or \emph{reshaping}) isomorphism corresponding to this partition as the linear map
\[
\bullet_{\{\mathcal{I}_1,\dots,\mathcal{I}_m\}} :
\mathrm{V}_1 \otimes \cdots \otimes \mathrm{V}_d \to 
\bigotimes_{k=1}^{m} \bigotimes_{j \in \mathcal{I}_k} \mathrm{V}_j,
\]
which rearranges and groups the vector spaces according to the partition, so that the codomain is viewed as a space of order-$m$ tensors. 
The action of this map on the rank-1 tensor $\tensor{B} = \bm{v}_1 \otimes \cdots \otimes \bm{v}_d \in \mathrm{V}_1 \otimes \dots \otimes \mathrm{V}_d$ is
\[
\bm{v}_1 \otimes \dots \otimes \bm{v}_d \mapsto \bigotimes_{k=1}^{m} \Bigl( \bigotimes_{j \in \mathcal{I}_k} \bm{v}_j \Bigr).
\]
Due to the universality of the tensor product~\cite{greub_multilinear_1978}, this fully defines the linear map $\bullet_{\{\mathcal{I}_1,\dots,\mathcal{I}_m\}}$.
The standard \emph{mode-$k$ flattening} corresponds to the partition $\mathcal{I}_1 = (k)$ and $\mathcal{I}_2 = (1,\dots,k-1,k+1,\dots,d)$.
If $m=2$, we call the flattening a \emph{matricization}.
For matricizations we also use the shorthand notation 
\(\tensor{A}_{\mathcal{I}_t} := \tensor{A}_{\{\mathcal{I}_t,\mathcal{I} \backslash \mathcal{I}_t\}}\).

Let $U_k : \mathrm{V}_k \to \mathrm{W}_k$ be linear maps between vector spaces $\mathrm{V}_k$ and $\mathrm{W}_k $ for $k=1,\dots,d$. 
The tensor product of linear maps, also called the Kronecker product, of $U_1, \dots, U_d$ is the linear map
\begin{align*}
	(U_1 \otimes \dots \otimes U_d)\colon \mathrm{V}_1 \otimes \dots \otimes \mathrm{V}_d &\to \mathrm{W}_1 \otimes \dots \otimes \mathrm{W}_d \quad \text{, such that}\\
(U_1 \otimes \dots \otimes U_d) (\bm{v}_1 \otimes \dots \otimes \bm{v}_d) &\;=\; (U_1\bm{v}_1) \otimes \dots \otimes (U_d\bm{v}_d).
\end{align*}
By universality \cite{greub_multilinear_1978}, this again fully defines the linear map $U_1 \otimes \dots \otimes U_d$. Representing $\tensor{A}$ and $U_k$ in coordinates, the map $\tensor{A} \mapsto (U_1 \otimes \dots \otimes U_d)(\tensor{A})$ is called a \emph{multilinear multiplication} and will be abbreviated to $(U_1, \dots,U_d) \cdot \tensor{A}$.
It follows from the foregoing that if $\mathcal{Y} = (U_1, \cdots, U_d) \cdot \tensor{A}$, then the mode-$k$ flattening satisfies
\begin{equation*}
	\label{eq:flattening_mult}
\tensor{Y}_{(k)} = U_k \tensor{A}_{(k)} (U_1 \otimes \dots \otimes U_{k-1} \otimes U_{k+1} \otimes \dots \otimes U_d)^\top
\end{equation*}
for each mode $k = 1, \dots,d$.
When all but one map $U_k$ are the identity $I$, the multilinear multiplication is shortened to $U_k \cdot_k \tensor{A} := (I,\dots, I, U_k,I,\dots,I)\cdot \tensor{A}$, where $U_k$ appears in the $k$th position. More generally, if $\mathcal{I}_k \subseteq (1,\dots,d)$ and $M \colon \bigotimes_{i \in \mathcal{I}_k} \mathrm{V}_i \to \bigotimes_{i \in \mathcal{I}_k} \mathrm{W}_i$, then $(M \cdot_{\mathcal{I}_k} \tensor{A})_{\mathcal{I}_k} := M \tensor{A}_{\mathcal{I}_k}$.
\subsection{Tucker decomposition}
\label{sec:tucker}
Let $\tensor{A} \in \mathrm{V}_1 \otimes \cdots \otimes \mathrm{V}_d$ be a tensor of order $d$, with $n_k$-dimensional vector spaces $\mathrm{V}_k$ for $k=1,\dots,d$. 
For each mode $k$, let $\mathrm{W}_k \subseteq \mathrm{V}_k$ be a subspace of dimension $r_k \leq n_k$, and denote by $U_k : \mathrm{W}_k \hookrightarrow \mathrm{V}_k$ the canonical inclusion.
An exact \emph{Tucker decomposition}~\cite{tucker_mathematical_1966} of $\tensor{A}$ is a representation of the form
\[
\tensor{A} = (U_1, \dots, U_d) \cdot \tensor{G},
\]
where $\tensor{A}$ is expressed as the image of a \emph{core tensor} $\tensor{G} \in \mathrm{W}_1 \otimes \cdots \otimes \mathrm{W}_d$ under the tensor product of inclusion maps. Note this simply states that $\tensor{A} \in \mathrm{W}_1 \otimes \dots \otimes \mathrm{W}_d$ and the Tucker decomposition records how $\tensor{A}$ is embedded into $\mathrm{V}_1 \otimes \dots \otimes \mathrm{V}_d$ via its inclusions $U_k$. The \emph{multilinear rank} of $\tensor{A}$ is the tuple $(r_1,\dots,r_d)$ where $r_k=\dim W_k$ and $\mathrm{W}_1 \otimes \dots \otimes \mathrm{W}_d$ is the smallest tensor product space containing $\tensor{A}$. The core tensor of the Tucker decomposition can be computed from $\tensor{A}$ by $\tensor{G} = ( U_1^{\dagger}, \dots, U_d^{\dagger} ) \cdot \tensor{A}$,
where $U_k^\dagger$ is a left inverse of $U_k$. If $\tensor{A}$ is represented in some coordinates with some tensor product basis of $\mathrm{V}_1 \otimes \dots \otimes \mathrm{V}_d$, then if $n = \max(n_1,\dots,n_d)$ and $r = \max(r_1,\dots,r_d)$, the storage requirement for the Tucker decomposition is bounded by \(\mathcal{O}(d nr + r^d)\).

In coordinates, the \emph{Higher-Order Singular Value Decomposition} (HOSVD)~\cite{de_lathauwer_multilinear_2000} is a specific Tucker decomposition method which chooses the left singular vectors of the mode-$k$ flattening $\tensor{A}_{(k)}$ as the columns of the matrix $U_k$.

A \emph{truncated HOSVD} (T-HOSVD) collects only the leading $r_k$ left singular vectors of $\tensor{A}_{(k)}$ in an $n_k \times r_k$ matrix $\widehat{U}_k$ to approximate $\tensor{A}$ by a lower multilinear rank tensor. 
That is, for a given multilinear rank, $(r_1,\dots,r_d)$, the \emph{T-HOSVD approximation} is
\[
\tensor{A} \approx 
(\widehat{U}_1 \widehat{U}_1^\top,\dots,\widehat{U}_d \widehat{U}_d^\top) \cdot \tensor{A}.
\]
In other words, $\tensor{A}$ is approximated by orthogonally projecting it to the tensor product subspace $\operatorname{span}(\widehat{U}_1 \otimes \dots \otimes \widehat{U}_d)$.
Computing the $\widehat{U}_k$ can be carried out in perfect parallelism. Consequently, this method will also be referred to as the \emph{parallel T-HOSVD}.

A \emph{sequentially truncated HOSVD} (ST-HOSVD) \cite{vannieuwenhoven_new_2012, hackbusch_tensor_2019} computes the subspaces $\mathrm{W}_k$ and the core tensor sequentially across the modes.
Each SVD is performed on the mode-$k$ flattening of the previously truncated tensor, leading to an $n_k \times r_k$ matrix $\overline{U}_k$ containing the $r_k$ left singular vectors. 
These sequential computations on previously truncated tensors can significantly reduce computational complexity and memory usage, particularly for large tensors. The resulting \emph{ST-HOSVD approximation} $(\overline{U}_1 \overline{U}_1^\top,\dots,\overline{U}_d \overline{U}_d^\top) \cdot \tensor{A}$ approximates $\tensor{A}$ by an orthogonal projection onto the tensor product subspace $\operatorname{span}(\overline{U}_1 \otimes \dots \otimes \overline{U}_d)$. 

We will refer to $(\breve{U}_1 \breve{U}_1^\top,\dots,\breve{U}_d \breve{U}_d^\top) \cdot \tensor{A}$ as an \emph{orthogonal Tucker projection} of $\tensor{A}$ if each $\breve{U}_i$ has orthogonal columns, e.g., both T-HOSVD and ST-HOSVD approximations are orthogonal Tucker projections.

In general, the subspaces computed by T-HOSVD and ST-HOSVD can be distinct.
Nonetheless, for a tensor $\tensor{A} \in \mathrm{V}_1 \otimes \cdots \otimes \mathrm{V}_d$, both produce \emph{quasi-optimal approximations}. This means these approximations are within some \emph{quasi-optimality factor} of the best approximation $\tensor{A}^{\mathrm{best}}$ with multilinear rank $(r_1,\dots,r_d)$~\cite{hackbusch_tensor_2019,grasedyck_hierarchical_2010}. For both T-HOSVD and ST-HOSVD approximations, $\tensor{A}^{\mathrm{S/T}}$, of $\tensor{A}$, it follows from the analysis in \cite{hackbusch_tensor_2019} that if truncation is performed only in the modes $\mathcal{I}_t \subseteq \mathcal{I}$, i.e., for all $k \notin \mathcal{I}_t$ the projection is onto the subspace $\mathrm{W}_k= \operatorname{span}(\tensor{A}_{(k)}) \subseteq \mathrm{V}_k$, then \begin{equation}
		\label{eq:quasi-optimal_HOSVD}
		\|\tensor{A} - \tensor{A}^\mathrm{S/T} \| \le \sqrt{|\mathcal{I}_t|} \|\tensor{A} - \tensor{A}^{\mathrm{best}}\|.
	\end{equation}
The value $1 \le \sqrt{|\mathcal{I}_t|} \le \sqrt{d}$ is the quasi-optimality factor of the T-HOSVD and ST-HOSVD approximations. Randomized SVDs have already been used to speed up both T-HOSVD and ST-HOSVD algorithms~\cite{minster_randomized_2020}.

\subsection{Mixed ST-HOSVD (M-HOSVD)}
With appropriate notational conventions, both T-HOSVD and ST-HOSVD can be treated uniformly. Let $\pi_k \colon \mathrm{V}_1 \otimes \cdots \otimes \mathrm{V}_d \to \mathrm{V}_1 \otimes \cdots \otimes \mathrm{V}_d$ denote the rank-$r$ truncated SVD projection that takes $\tensor{X} \mapsto U_k U_k^\top \cdot_k \tensor{X}$ where $U_k$ contains the $r$ dominant left singular vectors of $\tensor{X}_{(k)}$. Then, for $\mathcal{J} = (j_1,\dots,j_m) \subseteq (1,\dots,d)$, the ST-HOSVD projection is
\[
	P_\mathcal{J} := \pi_{j_m} \circ \dots \circ \pi_{j_1} \colon \mathrm{V}_1 \otimes \cdots \otimes \mathrm{V}_d \to \mathrm{V}_1 \otimes \cdots \otimes \mathrm{V}_d.
\]
For any fixed $\tensor{A} \in \mathrm{V}_1 \otimes \cdots \otimes \mathrm{V}_d$, \cref{sec:tucker} showed that there exists a linear orthogonal projector $\Pi^\tensor{A}_i$ such that $\pi_k(\tensor{A}) = \Pi^{\tensor{A}}_k(\tensor{A})$. Hence,
\[
	P_\mathcal{J}(\tensor{A}) = \Pi^{\tensor{A}}_\mathcal{J}(\tensor{A}) := \Pi^{\tensor{A}}_{j_m}\cdot_{j_m} \dots \Pi^{\tensor{A}}_{j_1}\cdot_{j_1} \tensor{A}.
\]
Given a partition $\mathcal{J}_1 \sqcup \dots \sqcup \mathcal{J}_n = (1,\dots,d)$, define 
\begin{align*}
	(P_{\mathcal{J}_1} \boxtimes \dots \boxtimes P_{\mathcal{J}_n})(\tensor{A}) &:= \Pi^{\tensor{A}}_{\mathcal{J}_1} \cdot_{\mathcal{J}_1} \dots \Pi^{\tensor{A}}_{\mathcal{J}_n} \cdot_{\mathcal{J}_n} \tensor{A}, \quad \text{where}\\
	\Pi^{\tensor{A}}_{\mathcal{J}_\ell} \cdot_{\mathcal{J}_\ell} \tensor{X} &:= \Pi^{\tensor{A}}_{j_{\ell,1}}\cdot_{j_{\ell,1}}\dots
	\Pi^{\tensor{A}}_{j_{\ell,m_\ell}}\cdot_{j_{\ell,m_\ell}} \tensor{X}
\end{align*} 
for each $\mathcal{J}_\ell = (j_{\ell,1}, \dots, j_{\ell,m_\ell})$, i.e., each projector $\Pi^{\tensor{A}}_{\mathcal{J}_\ell}$ is computed on $\tensor{A}$. For our later derivations, we define a mixed truncation strategy that combines both parallel and sequential truncation, mixing tensor products and compositions.
\begin{definition}[Mixed ST-HOSVD]
	\label{def:M-HOSVD}
	Let $\tensor{A}\in \mathrm{V}_1 \otimes \cdots \otimes \mathrm{V}_d$ and $\mathcal{J}_1 \sqcup \dots \sqcup \mathcal{J}_n := (1,\dots,d).$ 
	Let $\mathcal{J}_\ell = (j_1, \dots, j_m)$ and $P_{\mathcal{J}_\ell}$ denote the ST-HOSVD projection \[
	P_{\mathcal{J}_\ell}(\tensor{A}) := P_{\mathcal{J}_\ell} \cdot_{\mathcal{J}_\ell} \tensor{A} := (\pi_{j_m} \circ \dots \circ \pi_{j_1})(\tensor{A}),
	\] where each $\pi_{j_k}(\tensor{X})$ is a truncated SVD projection of $\tensor{X}_{(j_k)}$ to the desired rank.
	Then, the Mixed ST-HOSVD (M-HOSVD) projection is $(P_{\mathcal{J}_1} \boxtimes \dots \boxtimes P_{\mathcal{J}_n})(\tensor{A})$.
\end{definition} 
The M-HOSVD with a single group $\mathcal{J}_1 = \{1,\dots,d\}$ is equal to the ST-HOSVD, while the groupings $\mathcal{J}_\ell = \{\ell\}$ for $\ell = 1,\dots,d$ yields the T-HOSVD.
The quasi-optimality of ST-HOSVD and T-HOSVD implies quasi-optimality of the M-HOSVD as well.
\begin{lemma}
	\label{lemma:MTHOSVD_err}
	Let $\tensor{A} \in \mathrm{V}_1 \otimes \cdots \otimes \mathrm{V}_d$. Consider its multilinear rank-$(r_1, \dots, r_d)$ M-HOSVD approximation, $\widehat{\tensor{A}}$, using the partition $\bigsqcup_{\ell = 1}^n \mathcal{J}_{\ell} = (1,\dots,d)$ and let $m =  \bigl|\{\, r_k \mid r_k < \operatorname{rank}(\tensor{A}_{(k)}),\ k = 1,\dots,d \,\}\bigr|$. Then,
	\[
	\|\tensor{A} - \widehat{\tensor{A}}\| \le \sqrt{m} \|\tensor{A} - \tensor{A}^{\mathrm{best}}\|,
	\]
	where $\tensor{A}^{\mathrm{best}}$ is the best Tucker approximation of $\tensor{A}$ with the same multilinear rank.
	\begin{proof}
		Let $P_{\mathcal{J}_\ell}$ be the ST-HOSVD projection from \cref{def:M-HOSVD} for every $\ell = 1,\dots,n$ so that $P := P_{\mathcal{J}_1} \boxtimes \dots \boxtimes P_{\mathcal{J}_n}$. Then, by the usual telescoping argument as in \cite[Theorem 5.1]{vannieuwenhoven_new_2012}
		\[
			\|\tensor{A} - P (\tensor{A}) \|^2 = \sum_{\ell = 1}^{n} \|(P_{\mathcal{J}_1}, \dots, P_{\mathcal{J}_{\ell-1}}, P_{\mathcal{J}_{\ell}}^\perp, I, \dots, I)\cdot \tensor{A}\|^2 \le \sum_{\ell = 1}^{n}\|P_{\mathcal{J}_\ell}^\perp \cdot_{\mathcal{J}_\ell} \tensor{A}\|^2,
			\]
		where $P_{\mathcal{J}_\ell}^\perp:= I - P_{\mathcal{J}_\ell}$. The equality follows from the orthogonality of ST-HOSVDs across modes, and the inequality from the fact that the $P_{\mathcal{J}_\ell}(\tensor{A})$ project on disjoint modes and, hence, commute. Applying \cref{eq:quasi-optimal_HOSVD} to each summand concludes the proof.
	\end{proof}
\end{lemma}


\section{Hierarchical Tucker decomposition}
\label{sec:htd}
We now recall the definition and main properties of the HTD as found in~\cite{grasedyck_hierarchical_2010, hackbusch_new_2009, hackbusch_tensor_2019}.

\subsection{Definition}
\label{sec:htd-space}
Let $\tensor{A} \in \mathrm{V}_{1}\otimes \dots \otimes \mathrm{V}_{d}$. The tensor product is associative~\cite{greub_multilinear_1978}, so for any partition $\mathcal{I}_1 \sqcup \dots \sqcup \mathcal{I}_m = (1,\dots,d)$ we can regard $\tensor{A}$ as belonging to different
groupings of the vector spaces $\mathrm{V}_k$, via the flattening isomorphisms. For example, the partition $\mathcal{I}_1 = (1,2), \dots, \mathcal{I}_{m} = (d-1,d)$ permits viewing $\tensor{A}$ as the tensor of order $d/2$
\[
\tensor{A}_{\{\mathcal{I}_1,\dots,\mathcal{I}_{m}\}} \in (\mathrm{V}_{1}\otimes \mathrm{V}_{2})\otimes \cdots \otimes (\mathrm{V}_{d-1}\otimes \mathrm{V}_{d})
\]
 and $\mathcal{I}_1 = (1,\dots,d/2), \mathcal{I}_2 = (d/2+1,\dots,d)$ gives a view of $\tensor{A}$ as a second-order tensor
\[\tensor{A}_{\{\mathcal{I}_1,\mathcal{I}_2\}} \in (\mathrm{V}_{1} \otimes \cdots \otimes \mathrm{V}_{d/2})\otimes (\mathrm{V}_{d/2 + 1} \otimes \cdots \otimes \mathrm{V}_{d}).\]
Fundamentally, the HTD represents $\tensor{A}$ by
recursively grouping the vector spaces $\mathrm{V}_k$ according to a dimension tree, creating a hierarchy of nested subspaces.

\label{sec:htd-tree}
A \emph{dimension tree} $\mathcal{T}$ of $\mathcal{I} := (1,\dots,d)$ is a tree with $d$ leaves. Each node $t \in \mathcal{T}$ is identified by an integer starting from $1$ at the root, denoted by $\operatorname{root}(\mathcal{T})$, so that the number of any child is always larger than the parent. If the nodes are visited in ascending order, then this corresponds to a breadth-first traversal of the tree. The children of a node $t$ are denoted by $\operatorname{children}(t)$ and its parent by $\operatorname{parent}(t)$. 

To each node $t \in \mathcal{T}$, we attach two objects. First, an \emph{index set} $\mathcal{I}_t \subseteq \mathcal{I}$ satisfying:
	\begin{enumerate}
		\item if $t$ is a leaf, $\mathcal{I}_t$ is a singleton; 
		\item if $t$ is the root, $\mathcal{I}_t = \mathcal{I}$;
		\item if $t$ is an internal node, $\mathcal{I}_t$ is the disjoint union of the index sets of its children, i.e. $\mathcal{I}_t = \bigsqcup_{s \in \operatorname{children}(t)} \mathcal{I}_s$;
		\item $\mathcal{I}_t$ is sorted in ascending order.
	\end{enumerate}
Second, each $t$ has an associated tensor product space
\[
\mathrm{V}_{t}:=\bigotimes_{k\in\mathcal{I}_{t}}\mathrm{V}_{k}.
\]

The \emph{distance} between nodes in a dimension tree is the number of nodes minus one of the shortest path that connects these nodes.
A \emph{layer} or \emph{cut} in $\mathcal{T}$ is a subset of nodes $C \subset \mathcal{T}$ such that $\bigsqcup_{k\in C}\mathcal{I}_{k} = \mathcal{I}$. A tree can be \emph{sliced} into a tuple of layers so that every node is part of at least one layer. A layer can be identified with the leaves of a subtree rooted at $\operatorname{root}(\mathcal{T})$. The total number of layers in a slicing of $\mathcal{T}$ is denoted by $L$ and each layer $l = 1,\dots,L$ is denoted by $\mathcal{T}_l$. By definition, $\mathcal{T}_1$ and $\mathcal{T}_L$ contain the root and the leaves, respectively.
The matricization of $\tensor{A}$ in the modes corresponding to $t$ is abbreviated to $\tensor{A}_t := \tensor{A}_{\mathcal{I}_t}$ and the general flattening of $\tensor{A}$ using the index tuples of the $l$th layer is denoted by $\tensor{A}_{\mathcal{T}_l} := \tensor{A}_{\{\mathcal{I}_{t_1},\dots,\mathcal{I}_{t_m}\}}$, where $t_1, \dots, t_m$ are the $m$ nodes in $\mathcal{T}_l$. A possible slicing of a completely balanced tree $\mathcal{T}$ defines layers so that each $t \in \mathcal{T}_l$ is at a distance $l-1$ from $\operatorname{root}(\mathcal{T})$.

\label{sec:htd-rep}
An HTD of $\tensor{A} \in \mathrm{V}_1 \otimes \dots \otimes \mathrm{V}_d$ associates with each node $t \in \mathcal{T}$ a vector space
\begin{equation}
	\mathrm{W}_t \subseteq \bigotimes_{k\in \mathcal{I}_t} \mathrm{W}_k \subseteq \mathrm{V}_t
\end{equation}
of minimal dimension $r_t \le \operatorname{dim}(\mathrm{V}_t)$, also referred to as the \emph{hierarchical rank} at node $t$, such that $\tensor{A}_{\mathcal{T}_l} \in \mathrm{W}_{t_1} \otimes \dots \otimes \mathrm{W}_{t_m} \subseteq \mathrm{V}_{t_1} \otimes \dots \otimes \mathrm{V}_{t_m}$, where $\mathcal{T}_l = \{t_1, \dots, t_m\}$. This defines a hierarchy of Tucker decompositions of the flattenings of $\tensor{A}$ in successive layers of the dimension tree. That is, 
\[
\tensor{A} \in \bigcap_{l = 1}^L \bigotimes_{t\in \mathcal{T}_l} \mathrm{W}_t.
\]
For each node $t$, there exists a canonical inclusion map $F_t \colon \mathrm{W}_t \hookrightarrow \mathrm{V}_t$, which in coordinates can be represented by an $n_{t}\times r_{t}$ matrix, also denoted by $F_t$, where $n_t = \operatorname{dim}(\mathrm{V}_t) = \Pi_{i \in \mathcal{I}_t} \operatorname{dim}(\mathrm{V}_i)$ and $r_t = \operatorname{dim}(\mathrm{W}_t)$. If $F_t$ only has orthonormal columns, we will denote it by $U_t$.
Since, for each internal node $t$, with $ \operatorname{children}(t) = \{c_1, \dots, c_m\}$, $\mathrm{W}_t \subseteq \bigotimes_{i = 1}^m \mathrm{W}_{c_i}$, we can succinctly express a basis of $\mathrm{W}_t$ as linear combinations of tensor products of basis vectors of the spaces $\mathrm{W}_{c_i}$ for $i \in \{1,\dots ,m\}$. Hence,
\begin{equation}
		\label{eq:rec_def}
		F_t = (F_{c_1} \otimes \dots \otimes F_{c_m})(B_t)
	\end{equation}
for some $r_{c_1} \cdots r_{c_m} \times r_t$ matrix $B_t$. This $B_t$ represents the inclusion map $\mathrm{W}_t \hookrightarrow \mathrm{W}_{c_1} \otimes \dots \otimes \mathrm{W}_{c_m}$ in coordinates, and can equivalently be interpreted as the so-called \emph{transfer tensor}~\cite{grasedyck_hierarchical_2010} $\tensor{B}_t \in {W}_{c_1} \otimes \dots \otimes \mathrm{W}_{c_m} \otimes \mathrm{W}_t^*$, where $\mathrm{W}_t^{*}$ denotes the dual vector space of $\mathrm{W}_t$. For notational convenience, we set $B_t := F_t$ for all leaves $t \in \mathcal{T}_L$. Thus, an HTD can be represented by the tuple $(\mathcal{T}, \bm{B})$, where $\bm{B} = (B_t)_{t \in \mathcal{T}}$. We will refer to the mode that corresponds to the factor $\mathrm{W}_t^*$ of $\tensor{B}_t$ by $e_t$, i.e., $B_t = (\tensor{B}_t)_{\{\mathcal{I}_t, e_t\}} \in \mathbb{R}^{r_{c_1} \cdots r_{c_m} \times r_t}$.

\subsection{Storage complexity}
\label{sec:htd-stor}
One advantage of computing the HTD of $\tensor{A}$ when $\operatorname{dim}(\mathrm{W}_t) \ll \operatorname{dim}(\mathrm{V}_t)$, is that in coordinates, the transfer tensors admit a much more compact representation.
Let $r_t = \operatorname{dim}(\mathrm{W}_t)$, $R = \max_{t\in\mathcal{T}} r_t$, $r = \max_{t \in \mathcal{T}_L} r_t$, and $n = \max\{\operatorname{dim}(\mathrm{V}_1),\dots, \operatorname{dim}(\mathrm{V}_d)\}$. The storage cost of the $d-1$ transfer tensors and $d$ factor matrices of the HTD format can then be bounded by $\mathcal{O}((d-1) R^3 + dnr)$.
Note that $R$ can be much higher than $r$, being the rank bound of the Tucker decomposition, as the dimension of each $\mathrm{W}_t \subseteq \bigotimes_{c \in \operatorname{children}(t)} \mathrm{W}_c$ is bounded as $ r_t \le \Pi_{c \in \operatorname{children}(t)} r_c$.

\subsection{Root-normal form}
\label{sec:svd_htd_form}
Several HTD normal forms have been defined, allowing for decoupled, nodewise normalization and truncation \cite{Grasedyck2024}. We will use a variant of the root-normal form of \cite{Grasedyck2024} to simplify our derivations later. For further information about these normal forms and their properties we refer readers to \cite{Grasedyck2024}.
\begin{definition}[Essentially, {\cite[Definition 2]{Grasedyck2024}}]
	\label{def:svd_htd}
	An HTD $(\mathcal{T}, \bm{B})$ of a tensor $\tensor{A}$ is in SVD-HTD form at layer $l$, if for every node $t \in \mathcal{T}_l$, $F_t$, which is recursively defined by \cref{eq:rec_def} with $F_k = B_k$ for $k \in \mathcal{T}_L$, has the form
	\begin{displaymath}
		F_t = U_t \Sigma_t,
	\end{displaymath}
	where $U_t \Sigma_t V_t^\top$ is a compact SVD of $\tensor{A}_t$. An HTD is in SVD-HTD form if it is in SVD-HTD form at all layers of $\mathcal{T}$.
\end{definition}
We can use \cite[Algorithm 3]{Grasedyck2024} to transform a general HTD into SVD-HTD form. For clarity, we present a version of this algorithm in our notation as \cref{alg:SVD_HTD}.\footnote{Note that our notation for $U$ and $V$ is switched relative to~\cite{Grasedyck2024}.} 

\begin{algorithm}[!tb]
	\caption{SVD-HTD Transformation (Essentially, \cite[Algorithm 3]{Grasedyck2024})}
	\begin{algorithmic}[1]
			\label{alg:SVD_HTD}
			\REQUIRE HTD of tensor $\tensor{A}$: ($\mathcal{T}, \bm{B}$).\\
			\COMMENT{\emph{Stage 1: Leaves-to-root orthonormalization}}
			\STATE Orthonormalization of the HTD using~\cite[Algorithm 3]{grasedyck_hierarchical_2010}\\
			\COMMENT{\emph{Stage 2: Bring into SVD-HTD form}}
			\FOR{$l \gets 1$ to $L-1$}
			\FOR{$t \in \mathcal{T}_l \backslash \mathcal{T}_L$}
			\STATE $\{t_1, \dots ,t_m\} \gets \operatorname{children}(t)$\\
			Let $(\tensor{B}_t)_{\{\mathcal{I}_{t_1},\dots, \mathcal{I}_{t_m}, e_t\}} = (U_{t_1},\dots, U_{t_m}, U_{e_t})\cdot \tensor{G}$ be an HOSVD of $\tensor{B}_t \in \mathrm{W}_{t_1} \otimes \dots \otimes \mathrm{W}_{t_m} \otimes \mathrm{W}^*_t$, so that $(\tensor{B}_t)_{\mathcal{I}_{t_k}} = U_{t_k} \Sigma_{t_k} V_{t_k}^\top$ for all $k=1,\dots,m$.\\
			\STATE $(\tensor{B}_t)_{\{\mathcal{I}_{t_1},\dots, \mathcal{I}_{t_m}, e_t\}} \gets (\Sigma_{t_1}^{-1} ,\dots, \Sigma_{t_m}^{-1}, I) \cdot \tensor{G}$\\
			\STATE $B_{t_k} \gets B_{t_k}U_{t_k} \Sigma_{t_k}$ 		for $k = 1,\dots, m$\\
			\ENDFOR
			\ENDFOR	
	\end{algorithmic}
\end{algorithm}

\subsection{Existing algorithms}
\label{sec:htd-comp}
As explained in \cref{sec:tucker}, there are two dominant paradigms for truncation: \emph{parallel} and \emph{sequential}~\cite{hackbusch_tensor_2019}. In parallel truncation, all necessary projections are computed from the original tensor. In sequential truncation, the projections are computed and applied one after the other. In the case of the Tucker decomposition this corresponds to T-HOSVD and ST-HOSVD, respectively~\cite{de_lathauwer_multilinear_2000,vannieuwenhoven_new_2012,hackbusch_tensor_2019}. In the HTD case, this is the difference between Grasedyck's~\cite[Algorithm 1]{grasedyck_hierarchical_2010} and~\cite[Algorithm 2]{grasedyck_hierarchical_2010}, respectively.

The state-of-the-art algorithms currently found in software libraries~\cite{kressner_algorithm_2014,kressner_recompression_2017} for computing the HTD of a full dense tensor in coordinates, proceed in a sequential LTR direction along the dimension tree $\mathcal{T}$, i.e., they implement Grasedyck's~\cite[Algorithm 2]{grasedyck_hierarchical_2010}. 
At the leaf layer, they perform a T-HOSVD on the original tensor, $\tensor{A} \approx (U_1 \otimes \dots \otimes U_d) (\tensor{C})$, after which the core tensor, $\tensor{C}$, is reshaped based on the mode partition of $\mathcal{T}_{L-1}$ and another T-HOSVD approximation, $\tensor{C} \approx (B_{1,2} \otimes \dots \otimes B_{d-1,d}) (\tensor{C'})$, is performed. 
The combined 2-layer approximation is then
\[\tensor{A} \approx \underbrace{(U_1 \otimes \dots \otimes U_d)}_{\text{Leaf bases}} \underbrace{(B_{1,2} \otimes \dots \otimes B_{d-1,d})}_{\text{Transfer tensors}} (\tensor{C'}).\]
This act of reshaping and performing Tucker decompositions is recursively applied to the core tensor, layer by layer, until the root of the tree is reached. Formally, these LTR algorithms can be summarized as follows.
 
\begin{definition}[LTR HTD projection]
	\label{def:LTR_projection}
	Let \[
	P_l := \bullet_{\mathcal{T}_l}^{-1} \circ \operatorname{T-HOSVD} \circ \; \; \bullet_{\mathcal{T}_l}\colon \mathrm{V}_1 \otimes \dots \otimes \mathrm{V}_d \to \mathrm{V}_1 \otimes \dots \otimes \mathrm{V}_d,
	\] i.e., reshaping using the index tuples of the nodes on layer $\mathcal{T}_l$, followed by a T-HOSVD projection, and finally reshaping back to the original shape.
	The LTR algorithm of Grasedyck~\cite[Algorithm 2]{grasedyck_hierarchical_2010} mathematically is a sequence of T-HOSVD projections $P_l$ 
	applied sequentially layer by layer starting from the leaf layer of the tree and moving upwards to the root layer:
	\[
	\tensor{A}^{\mathrm{LTR}} = (P_{1} \circ \cdots \circ P_{L}) (\tensor{A}).
	\]
\end{definition}

If randomized SVDs are used to approximate the dominant left singular subspaces in each T-HOSVD projection, the dominant computational cost of an efficient implementation of this LTR algorithm on a balanced binary tree is bounded by
\begin{equation}
	\label{eq:Comp_LTR_THOSVD}
	\mathcal{O}\left(
	dN^{d}r
	+\sum_{i=3}^{\log_{2}d} 2^{i-1}R^{2^{i}+1}
	+R^{5}
	\right),
\end{equation}
with the same $N$, $R$ and $r$ as in \cref{sec:htd-stor}.
\begin{remark}
	It is not strictly necessary to use T-HOSVD at each layer. Any orthogonal Tucker approximation can be used to approximate $\tensor{A}$ or any of the core tensors $\tensor{C}$. This is evident from the proofs and theorems from Grasedyck~\cite{grasedyck_hierarchical_2010}, which only rely on orthogonality of the bases.
An obvious refinement is then to use a more recently developed Tucker decomposition, such as the ST-HOSVD, effectively performing sequential truncations not only across layers of the tree but also within each layer itself.
This modification gives the same error bound (see \cref{remark:error_orthogonality}), but can drastically reduce the cost of the intermediate SVDs,
in particular at the leaf layer where $d$ SVDs of the full tensor are
normally required.
Under this improved scheme the computational complexity in \cref{eq:Comp_LTR_THOSVD} reduces to
\begin{equation}
	\label{eq:Comp_LTR_STHOSVD}
	\mathcal{O}\left(
	\sum_{i=1}^{d} N^{d-i+1}r^{i}
	+\sum_{i=4}^{d}R^{i+1}
	\right).
\end{equation}
\end{remark}

Beyond the standard LTR algorithm, several newer approaches for computing or updating HTD representations have emerged.

For high-dimensional data where the full dense tensor is unavailable, the HTBB scheme~\cite{ryzhakov_black-box_2024} was introduced. It iteratively traverses the dimension tree searching for the maximum-volume index selection to build surrogate models from sample function evaluations. It uses cross approximation to construct an HTD approximation without ever forming the full tensor explicitly. An earlier approach in the same spirit is the black-box HTD method of Ballani and Grasedyck~\cite{ballani_black_2013}, which uses adaptive cross approximation to build its HTD approximation.

Focus has recently also shifted to streaming applications.
The BHT-l2r method~\cite{AKSOY2027115360} targets batches or streams of similar tensors, compressing them directly in the HTD format through an LTR procedure. This avoids recomputing the decomposition for each new tensor individually. A related development is the TTNN (and its sequential variant STTNN), which adapts randomized Nyström techniques for Tucker decompositions to tree tensor networks~\cite{bucci_randomized_2024}. These methods provide single-pass, streamable low-rank approximations by sampling and sketching each node, and traverse the tree in an RTL approach while performing sequential truncation inside each layer.

\section{An RTL truncation algorithm}
\label{sec:rtl}
We now present our randomized HTD algorithms for general dense tensors, focusing on the RTL direction and employing sequential truncation across and inside layers of the dimension tree. 

\subsection{Motivation}
\label{sec:motivation}
To illustrate the difference between the LTR and RTL directions, consider a tensor $\tensor{A} \in \mathbb{R}^{N \times N \times N \times N}.$ We want to approximate it by an HTD with a completely balanced binary tree and all hierarchical ranks equal to $R$.

A sequential decomposition algorithm in the LTR direction begins by performing an ST-HOSVD on the leaf layer. 
This leads to a core tensor of size $R \times R \times R\times R$ and has a cost of $\mathcal{O}(N^4R + N^3R^2 + N^2R^3 + NR^4)$ if rank-$R$ truncated randomized SVDs are used to approximate the dominant left singular subspaces. Then, the core tensor is matricized to an $R^2\times R^2$ matrix. Its ST-HOSVD approximation can be computed by using one truncated randomized SVD with a cost of $\mathcal{O}(R^5)$. In total, this LTR scheme has a complexity of $\mathcal{O}(N^4R + N^3R^2 + N^2R^3 + NR^4 + R^5)$.

\begin{figure}[tb]
	\centering
	\resizebox{\textwidth}{!}{
		\begin{tikzpicture}[thick, 
			every node/.style={scale=1}, 
			tensor/.style={draw, circle, minimum size = 0.5cm}, 
			bond/.style={-}, 
			leg/.style={-}, 
			bd/.style={font=\small\ttfamily},
			arrowstyle/.style={->, thick, shorten >=0.3cm, shorten <=0.3cm}]
			
			\newcommand{\vertspace}{-2.5cm}
			
			\node[tensor, fill=green!20] (T) at (0,0) {$\tensor{A}$};
			
			\foreach \i/\angle in {1/25, 2/155, 3/205, 4/335} {
				\draw[leg] (T) -- ++(\angle:1.2cm) node[pos=1.4, font=\small] {$n_{\i}$};
			}
			
			\node[tensor, fill=green!20] (F1) at ($(T) + (0,\vertspace)$) {};
			\node (split1F1) [left=0.5cm of F1] {};
			\node (split2F1) [right=0.5cm of F1] {};
			\draw[leg] (F1) -- (split1F1.center) {};
			\draw[leg] (F1) -- (split2F1.center) {};
			\foreach \i/\angle in {1/135, 2/225} {
				\draw[leg] (split1F1.center) -- ++(\angle:0.6cm) node[pos=1.4, font=\small] {$n_{\i}$};
			}
			\foreach \i/\angle in { 3/45, 4/-45} {
				\draw[leg] (split2F1.center) -- ++(\angle:0.6cm) node[pos=1.4, font=\small] {$n_{\i}$};
			}

			\draw[arrowstyle] (T) -- (F1) node[midway, left,xshift=-0.5cm,align=left] {Reshape into \\ $(n_{1}n_{2} \times n_{3} n_{4})$};
			
			\node[tensor, fill=green!20] (F1S) at ($(F1) + (0,\vertspace)$) {$\Sigma$};
			\node[tensor, fill=red!20] (F1U) at ($(F1S) + (-1.5,0)$) {$U$};
			\node[tensor, fill=red!20] (F1V) at ($(F1S) + (1.5,0)$) {$V$};
			\node (splitUF1) [left=0.1cm of F1U] {};
			\node (splitVF1) [right=0.1cm of F1V] {};
			\draw[leg] (F1U) -- (F1S) node[midway,above] {$r_{34}$ };
			\draw[leg] (F1V) -- (F1S) node[midway,above] {$r_{12}$};
			\draw[leg] (F1U) -- (splitUF1.center) {};
			\draw[leg] (F1V) -- (splitVF1.center) {};
			\foreach \i/\angle in {1/135, 2/225} {
				\draw[leg] (splitUF1.center) -- ++(\angle:0.6cm) node[pos=1.4, font=\small] {$n_{\i}$};
			}
			\foreach \i/\angle in { 3/45, 4/-45} {
				\draw[leg] (splitVF1.center) -- ++(\angle:0.6cm) node[pos=1.4, font=\small] {$n_{\i}$};
			}
			
			\draw[arrowstyle] (F1) -- (F1S) node[midway, left,align=left] {SVD};

			\node[tensor, fill=green!20] (F1US) at ($(F1S) + (-4,\vertspace)$) {$\widehat{U} \widehat{\Sigma}$};
			\node (splitUSF1) [left=0.5cm of F1US] {};
			\node (splitotherF1) [right=0.5cm of F1US] {};
			\draw[leg] (F1US) -- (splitUSF1.center) {};
			\draw[leg] (F1US) -- (splitotherF1.center) node[midway,above] {$r_{12}$};
			\foreach \i/\angle in {1/135, 2/225} {
				\draw[leg] (splitUSF1.center) -- ++(\angle:0.6cm) node[pos=1.4, font=\small] {$n_{\i}$};
			}
			
			\node[tensor, fill=green!20] (F2US) at ($(F1S) + (4,\vertspace)$) {$\widehat{V} \widehat{\Sigma}$};
			\node (splitUSF2) [left=0.5cm of F2US] {};
			\node (splitotherF1T) [right=0.5cm of F2US] {};
			\draw[leg] (F2US) -- (splitUSF2.center) {};
			\draw[leg] (F2US) -- (splitotherF1T.center) node[midway,above] {$r_{34}$};
			\foreach \i/\angle in {3/135, 4/225} {
				\draw[leg] (splitUSF2.center) -- ++(\angle:0.6cm) node[pos=1.4, font=\small] {$n_{\i}$};
			}
			
			\draw[arrowstyle] (F1S) -- (F1US) node[midway, left,align=center] {Continue with\\ $\widehat{U} \widehat{\Sigma}$};
			\draw[arrowstyle] (F1S) -- (F2US) node[midway, right,align=center] {Continue with\\ $\widehat{V} \widehat{\Sigma}$};
			
			\node[tensor, fill=green!20] (F3) at ($(F1US) + (-2,\vertspace)$) {};
			\node (split1F3) [left=0.3cm of F3] {};
			\node (split2F3) [right=0.3cm of F3] {};
			\draw[leg] (F3) -- (split2F3.center) {};
			\draw[leg] (F3) -- (split1F3.center) node[pos=1.4, font=\small] {$n_{1}$};
			\draw[leg] (split2F3.center) -- ++(45:0.6cm) node[pos=1.4, font=\small] {$n_{2}$};
			\draw[leg] (split2F3.center) -- ++(-45:0.6cm) node[pos=1.4, font=\small] {$r_{12}$};

			\node[tensor, fill=green!20] (F4) at ($(F1US) + (2,\vertspace)$) {};
			\node (split1F4) [left=0.3cm of F4] {};
			\node (split2F4) [right=0.3cm of F4] {};
			\draw[leg] (F4) -- (split2F4.center) {};
			\draw[leg] (F4) -- (split1F4.center) node[pos=1.4, font=\small] {$n_{2}$};
			\draw[leg] (split2F4.center) -- ++(45:0.6cm) node[pos=1.4, font=\small] {$n_{1}$};
			\draw[leg] (split2F4.center) -- ++(-45:0.6cm) node[pos=1.4, font=\small,xshift=0.3cm] {$r_{12}$};
			
			\draw[arrowstyle] (F1US) -- (F3) node[midway, left,xshift=-0.5cm,align=left] {Reshape into \\ $(n_{1} \times n_{2}r_{34})$};
			\draw[arrowstyle] (F1US) -- (F4) node[midway,xshift=0.3cm, right,align=right] {Reshape into \\ $(n_{2} \times r_{34}n_{1})$};
			
			\node[tensor, fill=green!20] (F5) at ($(F2US) + (2,\vertspace)$) {};
			\node (split1F5) [left=0.3cm of F5] {};
			\node (split2F5) [right=0.3cm of F5] {};
			\draw[leg] (F5) -- (split2F5.center) {};
			\draw[leg] (F5) -- (split1F5.center) node[pos=1.4, font=\small] {$n_{4}$};
			\draw[leg] (split2F5.center) -- ++(45:0.6cm) node[pos=1.4, font=\small] {$n_{3}$};
			\draw[leg] (split2F5.center) -- ++(-45:0.6cm) node[pos=1.4, font=\small] {$r_{34}$};

			\node[tensor, fill=green!20] (F6) at ($(F2US) + (-2,\vertspace)$) {};
			\node (split1F6) [left=0.3cm of F6] {};
			\node (split2F6) [right=0.3cm of F6] {};
			\draw[leg] (F6) -- (split2F6.center) {};
			\draw[leg] (F6) -- (split1F6.center) node[pos=1.4, font=\small] {$n_{3}$};
			\draw[leg] (split2F6.center) -- ++(45:0.6cm) node[pos=1.4, font=\small] {$n_{4}$};
			\draw[leg] (split2F6.center) -- ++(-45:0.6cm) node[pos=1.4, font=\small] {$r_{34}$};
			
			\draw[arrowstyle] (F2US) -- (F5) node[midway, right,xshift=0.5cm,align=right] {Reshape into \\ $(n_{4} \times r_{12}n_{3})$};
			\draw[arrowstyle] (F2US) -- (F6) node[midway, left,xshift=-0.3cm,align=left] {Reshape into \\ $(n_{3} \times n_{4}r_{12})$};
			
			\node[tensor, fill=green!20] (F3S) at ($(F3) + (0,\vertspace)$) {$\Sigma$};
			\node[tensor, fill=red!20] (F3U) at ($(F3S) + (-1,0)$) {$U$};
			\node[tensor, fill=red!20] (F3V) at ($(F3S) + (1,0)$) {$V$};
			\node (splitUF3) [left=0.1cm of F3U] {};
			\node (splitVF3) [right=0.1cm of F3V] {};
			\draw[leg] (F3U) -- (F3S) node[midway,above] {$r_{1}$ };
			\draw[leg] (F3V) -- (F3S) node[midway,above] {$r_{1}$};
			\draw[leg] (F3U) -- (splitUF3.center) node[pos=1.7, font=\small] {$n_{1}$};
			\draw[leg] (F3V) -- (splitVF3.center) {};
			\draw[leg] (splitVF3.center) -- ++(45:0.4cm) node[pos=1.6, font=\small] {$n_{2}$};
			\draw[leg] (splitVF3.center) -- ++(-45:0.4cm) node[pos=1.6, font=\small] {$r_{12}$};
			
			\node[tensor, fill=green!20] (F4S) at ($(F4) + (0,\vertspace)$) {$\Sigma$};
			\node[tensor, fill=red!20] (F4U) at ($(F4S) + (-1,0)$) {$U$};
			\node[tensor, fill=red!20] (F4V) at ($(F4S) + (1,0)$) {$V$};
			\node (splitUF4) [left=0.1cm of F4U] {};
			\node (splitVF4) [right=0.1cm of F4V] {};
			\draw[leg] (F4U) -- (F4S) node[midway,above] {$r_{2}$ };
			\draw[leg] (F4V) -- (F4S) node[midway,above] {$r_{2}$};
			\draw[leg] (F4U) -- (splitUF4.center) node[pos=1.7, font=\small] {$n_{2}$};
			\draw[leg] (F4V) -- (splitVF4.center) {};
			\draw[leg] (splitVF4.center) -- ++(45:0.4cm) node[pos=1.6, font=\small] {$n_{1}$};
			\draw[leg] (splitVF4.center) -- ++(-45:0.4cm) node[pos=1.6, font=\small] {$r_{12}$};
			
			\node[tensor, fill=green!20] (F5S) at ($(F5) + (0,\vertspace)$) {$\Sigma$};
			\node[tensor, fill=red!20] (F5U) at ($(F5S) + (-1,0)$) {$U$};
			\node[tensor, fill=red!20] (F5V) at ($(F5S) + (1,0)$) {$V$};
			\node (splitUF5) [left=0.1cm of F5U] {};
			\node (splitVF5) [right=0.1cm of F5V] {};
			\draw[leg] (F5U) -- (F5S) node[midway,above] {$r_{4}$ };
			\draw[leg] (F5V) -- (F5S) node[midway,above] {$r_{4}$};
			\draw[leg] (F5U) -- (splitUF5.center) node[pos=1.7, font=\small] {$n_{4}$};
			\draw[leg] (F5V) -- (splitVF5.center) {};
			\draw[leg] (splitVF5.center) -- ++(45:0.4cm) node[pos=1.6, font=\small] {$n_{3}$};
			\draw[leg] (splitVF5.center) -- ++(-45:0.4cm) node[pos=1.6, font=\small] {$r_{34}$};
			
			\node[tensor, fill=green!20] (F6S) at ($(F6) + (0,\vertspace)$) {$\Sigma$};
			\node[tensor, fill=red!20] (F6U) at ($(F6S) + (-1,0)$) {$U$};
			\node[tensor, fill=red!20] (F6V) at ($(F6S) + (1,0)$) {$V$};
			\node (splitUF6) [left=0.1cm of F6U] {};
			\node (splitVF6) [right=0.1cm of F6V] {};
			\draw[leg] (F6U) -- (F6S) node[midway,above] {$r_{3}$ };
			\draw[leg] (F6V) -- (F6S) node[midway,above] {$r_{3}$};
			\draw[leg] (F6U) -- (splitUF6.center) node[pos=1.7, font=\small] {$n_{3}$};
			\draw[leg] (F6V) -- (splitVF6.center) {};
			\draw[leg] (splitVF6.center) -- ++(45:0.4cm) node[pos=1.6, font=\small] {$n_{4}$};
			\draw[leg] (splitVF6.center) -- ++(-45:0.4cm) node[pos=1.6, font=\small] {$r_{34}$};
			
			\draw[arrowstyle] (F3) -- (F3S) node[midway, left,align=left] {SVD};
			\draw[arrowstyle] (F4) -- (F4S) node[midway, right,align=right] {SVD};	
			\draw[arrowstyle] (F5) -- (F5S) node[midway, right,align=right] {SVD};
			\draw[arrowstyle] (F6) -- (F6S) node[midway, left,align=left] {SVD};
			
			\node[tensor, fill=green!20] (F3US) at ($(F3S) + (0,\vertspace)$) {$\widehat{U}\widehat{\Sigma}$};
			\node (U1l) [left=0.4cm of F3US] {};
			\node (U1r) [right=0.4cm of F3US] {};
			\draw[leg] (F3US) -- (U1r) node[midway,above] {$r_{1}$ };
			\draw[leg] (F3US) -- (U1l) node[midway,above] {$n_{1}$ };
			
			\node[tensor, fill=green!20] (F4US) at ($(F4S) + (0,\vertspace)$) {$\widehat{U}\widehat{\Sigma}$};
			\node (U2l) [left=0.4cm of F4US] {};
			\node (U2r) [right=0.4cm of F4US] {};
			\draw[leg] (F4US) -- (U2r) node[midway,above] {$r_{2}$ };
			\draw[leg] (F4US) -- (U2l) node[midway,above] {$n_{2}$ };
			
			\node[tensor, fill=green!20] (F5US) at ($(F5S) + (0,\vertspace)$) {$\widehat{U}\widehat{\Sigma}$};
			\node (U3l) [left=0.4cm of F5US] {};
			\node (U3r) [right=0.4cm of F5US] {};
			\draw[leg] (F5US) -- (U3r) node[midway,above] {$r_{4}$ };
			\draw[leg] (F5US) -- (U3l) node[midway,above] {$n_{4}$ };

			\node[tensor, fill=green!20] (F6US) at ($(F6S) + (0,\vertspace)$) {$\widehat{U}\widehat{\Sigma}$};
			\node (U4l) [left=0.4cm of F6US] {};
			\node (U4r) [right=0.4cm of F6US] {};
			\draw[leg] (F6US) -- (U4r) node[midway,above] {$r_{3}$ };
			\draw[leg] (F6US) -- (U4l) node[midway,above] {$n_{3}$ };
			
			\draw[arrowstyle] (F3S) -- (F3US) node[midway, left,align=left] {};
			\draw[arrowstyle] (F4S) -- (F4US) node[midway, right,align=right] {};
			\draw[arrowstyle] (F5S) -- (F5US) node[midway, right,align=right] {};
			\draw[arrowstyle] (F6S) -- (F6US) node[midway, left,align=left] {};
		\end{tikzpicture}
	}
\footnotesize
	\caption{\footnotesize Diagram showing the RTL algorithm for an order-4 tensor $\tensor{A}$. Red nodes have orthonormal columns, while green nodes contain the (approximate) singular values of $\tensor{A}$ at that layer.}
	\label{fig:RTL_algorithm}
\end{figure}
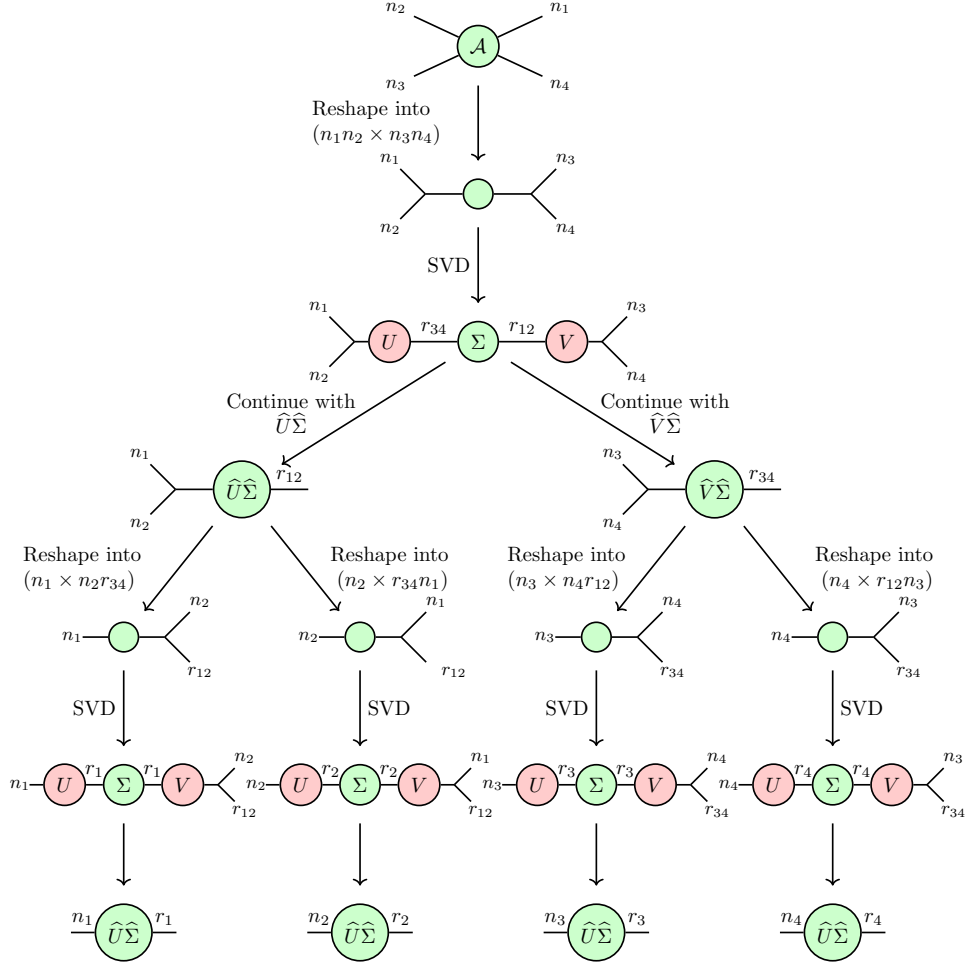

By contrast, a simplified version of the proposed sequential RTL truncation algorithm, which is described in \cref{sec:algorithm}, proceeds as shown in \cref{fig:RTL_algorithm} for this concrete case. First, a rank-$R$ truncated randomized SVD is performed on the $N^2\times N^2$ matricization of $\tensor{A}$, with a cost of $\mathcal{O}(N^4R)$. The left and right singular vectors are then each reshaped into their two possible $N \times NR$ matricizations. We then compute the dominant left singular subspaces of each of these $4$ matricizations using rank-$R$ truncated randomized SVDs, giving a cost of $\mathcal{O}(4N^2R^2)$. The proposed RTL algorithm thus has a total cost of $\mathcal{O}(N^4R + 4N^2R^2)$. Both algorithms share the same dominant cost $\mathcal{O}(N^{4}R)$, but their secondary terms differ substantially. RTL replaces a factor $N^{2}$ by $R$, whereas LTR replaces only one factor $N$. This gap widens for higher-order tensors, giving RTL a lower complexity in this uniform rank case. 

\subsection{A mathematical algorithm}
\label{sec:meth}
The idea of RTL algorithms is conceptually straightforward: perform Tucker decompositions sequentially in the RTL direction instead of LTR.
\begin{definition}[RTL HTD projection]\label{def:RTL_projection}
	Let \[
P_l := \bullet_{\mathcal{T}_l}^{-1} \circ \operatorname{Tucker} \circ \; \; \bullet_{\mathcal{T}_l}\colon \mathrm{V}_1 \otimes \dots \otimes \mathrm{V}_d \to \mathrm{V}_1 \otimes \dots \otimes \mathrm{V}_d,
\] be a projection that consists of reshaping using the index tuples of the nodes on layer $\mathcal{T}_l$, followed by an orthogonal Tucker projection with quasi-optimality factor $C_l$, and finally reshaping back to the original shape. An RTL approximation is a sequence of projections $P_l$ 
	applied layer by layer starting from the root layer of the tree and moving downwards to the leaf layer:
	\[
	\tensor{A}^{\mathrm{RTL}} = (P_{L} \circ \cdots \circ P_{1}) (\tensor{A}).
	\]
\end{definition}

\begin{remark}
	The projection scheme above is different from the parallel ``RTL'' algorithm introduced by Grasedyck \cite[Algorithm 1]{grasedyck_hierarchical_2010}. Our definition employs \emph{sequential} projections, which are computed from the previously projected tensor. The algorithm described by Grasedyck computes the projections from the original tensor. As a result, the order in which nodes are visited is inconsequential, enabling \emph{parallel} calculations. In this paper, an RTL algorithm is defined as one that enforces a strict traversal order, beginning at the root and progressing toward the leaf nodes and vice versa for LTR. As such, only algorithms following the sequential paradigm can thus be classified as either LTR or RTL. If orthogonal Tucker projections are used, however,\cite[Lemma 3.15]{grasedyck_hierarchical_2010} still ensures the ranks of the previous nodes do not increase after successive projections, leading to an HTD. 
\end{remark}

\begin{remark}
	\cref{def:RTL_projection} is essentially the sequential RTL HTD truncation scheme of Hackbusch \cite[(11.58)]{hackbusch_tensor_2019}. However, no numerical results were reported and no randomized SVDs were used in~\cite{hackbusch_tensor_2019}, and the method was seemingly abandoned in favor of LTR~\cite{grasedyck_hierarchical_2010}, likely due to the unfavorable computational complexity of RTL in absence of randomized SVDs.
\end{remark}

\subsection{Truncation and error analysis}
\label{sec:err}
In exact arithmetic and without truncation, the order in which projections are applied to obtain an HTD representation, whether in the LTR or RTL direction, i.e., \cref{def:LTR_projection} and \cref{def:RTL_projection}, respectively, does not affect the final projected tensor. Indeed, in the absence of truncation projections act as identities: $P_i(\tensor{A}) = \tensor{A}$. Hence, the order of ``projecting'' is irrelevant.
This property is analogous to the equivalence between ST-HOSVD and T-HOSVD in the absence of truncation. 
However, this equivalence breaks down with truncation.

In the presence of truncation, the LTR direction is definitely the more intuitive choice because it preserves the nestedness of the subspaces $\mathrm{W}_t$ as shown in \cite[p. 440]{hackbusch_tensor_2019}. By contrast, RTL projections do not guarantee nestedness once truncation is introduced, hence creating an additional source of approximation errors.

\begin{theorem}
	\label{thm:bad_RTL_error}
	Let $\tensor{A} \in \mathrm{V}_1 \otimes \dots \otimes \mathrm{V}_d$. The sequential RTL algorithm from \cref{def:RTL_projection} applied to $\tensor{A}$, where each orthogonal Tucker projector $P_l$ has quasi-optimality factor $C_l$, satisfies the following error bound:
	\begin{equation*}
		\label{eq:RTL_error_bound}
		\|\tensor{A} - (P_L \circ \cdots \circ P_1)(\tensor{A})\|
		\le \sum_{l=1}^{L} \|(P_l^\perp \circ P_{l-1} \circ \cdots \circ P_1) (\tensor{A})\|
		\le  \|\tensor{A} - \tensor{A}^{\mathrm{best}}\| \sum_{l=1}^{L} \tilde{C}_{l},
	\end{equation*}
	where $\tilde{C}_{l}$ is a computable constant for each layer based on $C_k$ for $k \le l$, $\tensor{A}^{\mathrm{best}}$ is the best HTD approximation of $\tensor{A}$ with the same hierarchical rank as $(P_L \circ \cdots \circ P_1)(\tensor{A})$, and $P_l^\perp(\tensor{X}) := \tensor{X} - P_l(\tensor{X})$.

\begin{proof}
	Consider a simple projection chain $P_{i,j}:= (P_i \circ P_j) (\tensor{A})$. Then,
	\begin{align*}
		\|\tensor{A} - P_{i,j} (\tensor{A})\|
		&= \|\tensor{A} - P_j (\tensor{A} )+ P_j (\tensor{A}) - P_{i,j} (\tensor{A})\| \le \|P_j^\perp (\tensor{A})\| + \|(P_i^\perp \circ P_j) (\tensor{A})\|.
	\end{align*}
	Recursively applying this bound to $(P_L\circ \cdots \circ P_1) (\tensor{A})$ yields the first inequality. This corresponds to using the triangle inequality on a telescoping sum.
	
	We prove on induction that $\|(P_l^\perp \circ P_{l-1} \circ \cdots \circ P_1 )(\tensor{A})\| \le \tilde{C}_{l} \|\tensor{A} - \tensor{A}^{\mathrm{best}}\|$ for each $l$. 
	In the base case $l=1$, $P_1$ is the identity projection because $\mathcal{T}_1 = \operatorname{root}(\mathcal{T})$ and $\mathcal{I}_{1} = (1,\dots,d)$, hence the bound trivially holds with $\tilde{C}_1 = 0$. 
	
	We assume the statement holds for layer $l$ and show that it also holds for $l+1$. The hierarchical ranks associated with the nodes on layer $l$ determine the multilinear rank of $\tensor{A}_{\mathcal{T}_l}$. We denote by $\mathcal{M}_{\mathcal{T}_l}$ the space of all tensors whose $\mathcal{T}_l$-flattening has multilinear rank bounded by these ranks. Furthermore, $\mathcal{H}_{\mathcal{T}_l}$ denotes the space of all tensors with bounded hierarchical rank $r_t$ for all nodes $t$ in the subtree consisting of the first $l$ layers, i.e., $\bigcup_{i=1}^l \mathcal{T}_i$, so that $\tensor{A}^{\mathrm{best}} \in \mathcal{H}_{\mathcal{T}_L}$. It is well known that $\mathcal{H}_{\mathcal{T}_L} \subseteq \mathcal{H}_{\mathcal{T}_{l+1}}\subset \mathcal{M}_{\mathcal{T}_{l+1}}$. Due to this nesting and the quasi-optimality of $P_{l+1}$, the following holds for $l+1$: 
	\[\|(P_{l+1}^\perp \circ P^{l})(\tensor{A})\| \le C_{l+1} \min_{\tensor{X} \in \mathcal{M}_{\mathcal{T}_{l+1}}} \|(P^{l} (\tensor{A}) -\tensor{X})\| \le C_{l+1} \min_{\tensor{X} \in \mathcal{H}_{\mathcal{T}_{L}}} \|P^{l}(\tensor{A})-\tensor{A} + \tensor{A} - \tensor{X}\|,\]
	where $P^{l} := P_l\circ \cdots \circ P_1$. Using a telescoping sum of $P^{i}(\tensor{A})$ terms with decreasing $i$ and applying the triangle inequality, we can proceed as follows
	\begin{align*}
		\|(P_{l+1}^\perp \circ P^l) (\tensor{A})\| &\le C_{l+1}\left(\sum_{i=2}^{l} \|(P_{i}^\perp \circ P^{i-1})(\tensor{A})\| + \min_{\tensor{X} \in \mathcal{H}_{\mathcal{T}_{L}}}\|\tensor{A}-\tensor{X}\|\right) \\
		&\le \underbrace{C_{l+1}\left(1 + \sum_{i=2}^{l} \tilde{C}_{i}\right)}_{\tilde{C}_{l+1}} \|\tensor{A} - \tensor{A}^{\mathrm{best}}\|,
	\end{align*}
where the last step used the induction hypothesis. This concludes the proof.
\end{proof}
\end{theorem}
\begin{remark}
	The $\tilde{C}_{l}$ from \cref{thm:bad_RTL_error} can be computed recursively starting from layer $2$, using the quasi-optimality factor of each $P_l$, $\tilde{C}_1 = 0$ and the formula $\tilde{C}_{l} = C_{l}(1 + \sum_{i=1}^{l-1} \tilde{C}_{i})$. These $\tilde{C}_{l}$ do, however, increase exponentially in the number of layers. For example, for a completely balanced dimension tree with $d=8$ ($L = 4$) and where each $C_l = \sqrt{|\mathcal{T}_l|}$, the $\tilde{C}_{l}$ would, in order, be $[0, \sqrt{2}, 2 + 2\sqrt{2}, 6(2+ \sqrt{2}) ]$.
\end{remark}

Fortunately, if an SVD-based Tucker decomposition is used for the quasi-optimal projector $P_l$ at each layer, then we have $\|(P_l^\perp \circ P_{l-1} \circ \cdots \circ P_1) (\tensor{A})\| \le \|P_l^\perp(\tensor{A})\|$, which is a consequence of the next technical result. Although the core result appears in~\cite[Theorem 11.63]{hackbusch_tensor_2019}, for self‑containedness, we include here a simpler proof.
\begin{lemma}
	\label{lemma:interlacing_singval}
	Let $\tensor{A} \in \mathrm{V}_1 \otimes \dots \otimes \mathrm{V}_d$ be a tensor. Let $\mathcal{I}_t \subseteq (1,\dots,d)$ and let $\Pi_{\mathcal{I}_t} = \widehat{U}_{\mathcal{I}_t} \widehat{U}_{\mathcal{I}_t}^\top$ be a rank-$r_t$ truncated SVD projection of the flattening $\tensor{A}_{\mathcal{I}_t}$. If $\widehat{\tensor{A}} := (\widehat{U}_{\mathcal{I}_t} \widehat{U}_{\mathcal{I}_t}^\top) \cdot_{\mathcal{I}_t} \tensor{A}$,
	then, for all index sets $\mathcal{I}_c$ where $ \mathcal{I}_c \subseteq \mathcal{I}_{t}$, $\mathcal{I}_t \subseteq \mathcal{I}_{c}$ or $\mathcal{I}_c \cap \mathcal{I}_{t} = \emptyset$, we have the elementwise inequality
	\[
		\sigma_{\mathcal{I}_c}(\widehat{\tensor{A}}) \le \sigma_{\mathcal{I}_c}(\tensor{A}),
	\]
	where $\sigma_{\mathcal{I}_c}(\tensor{A})$ denotes the vector of singular values of $\tensor{A}_{\mathcal{I}_c}$ sorted decreasingly.
	\begin{proof}
		 \textbf{Case} $\bm{\mathcal{I}_c \subseteq \mathcal{I}_{t}}$: Without loss of generality, we can assume that $\mathcal{I}_t = (1,\dots,k)$. Let $\bm{v} \in \bigotimes_{i \in \mathcal{I}_c} \mathrm{V}_i$. Then, $\|\bm{v}^\top \widehat{\tensor{A}}_{\mathcal{I}_c}\| = \|(\bm{v} \otimes I)^\top  P^\top \widehat{\tensor{A}}_{\mathcal{I}_t}\|$, with suitably sized identity matrix $I$ and a permutation matrix $P^\top$ that maps the tensor factors in $\mathcal{I}_c$ to $(1,\dots,l)$. We know by definition of $\widehat{\tensor{A}}$ that 
		\[\widehat{\tensor{A}}_{\mathcal{I}_t} = \widehat{U}_{\mathcal{I}_t} \widehat{\Sigma}_{\mathcal{I}_t} \widehat{V}_{\mathcal{I}_t}^\top = U_{\mathcal{I}_t}
		\begin{bmatrix}\widehat{\Sigma}_{\mathcal{I}_t} & \\ & 0\end{bmatrix}
		V_{\mathcal{I}_t}^\top,\] 
	where $\widehat{\Sigma}_{\mathcal{I}_t}$ are the singular values corresponding to $\widehat{U}_{\mathcal{I}_t}$ and $\widehat{V}_{\mathcal{I}_t}$ are the corresponding right singular vectors. Therefore, we have 
		
			\[\|\bm{v}^\top \widehat{\tensor{A}}_{\mathcal{I}_c}\| = \left\|(\bm{v} \otimes I)^\top P^\top U_{\mathcal{I}_t}
			\begin{bmatrix}\widehat{\Sigma}_{\mathcal{I}_t} & \\ & 0\end{bmatrix}
			V_{\mathcal{I}_t}^\top \right\|.\]
			From the definition of the Frobenius norm and its orthogonal invariance, we have that
			\[\left\|(\bm{v} \otimes I)^\top P^\top U_{\mathcal{I}_t}
			\begin{bmatrix}\widehat{\Sigma}_{\mathcal{I}_t} & \\ & 0\end{bmatrix}
			V_{\mathcal{I}_t}^\top \right\| \le \|(\bm{v} \otimes I)^\top P^\top U_{\mathcal{I}_t}\Sigma_{\mathcal{I}_t}V_{\mathcal{I}_t}^\top\| = \|(\bm{v} \otimes I)^\top P^\top \tensor{A}_{\mathcal{I}_t}\|,\]
		where $\Sigma_{\mathcal{I}_t}$ are all the singular values of $\tensor{A}_{\mathcal{I}_t}$. As this holds for all vectors $\bm{v}$, the proof is concluded by exploiting the max-min characterization of the $k$th singular value $\sigma_k$:
		\begin{align*}
			\sigma_k(\widehat{\tensor{A}}_{\mathcal{I}_c}) &= \max_{\operatorname{dim}(\mathrm{V}) =k} \min_{\substack{\bm{v} \in \mathrm{V} \\ \|\bm{v}\| = 1}} 	\|\bm{v}^\top \widehat{\tensor{A}}_{\mathcal{I}_c}\|
			\le 	\max_{\operatorname{dim}(\mathrm{V}) =k} \min_{\substack{\bm{v} \in \mathrm{V} \\ \|\bm{v}\| = 1}} 	\|\bm{v}^\top \tensor{A}_{\mathcal{I}_c}\|
			= \sigma_k (\tensor{A}_{\mathcal{I}_c}).
		\end{align*}
	
		\textbf{Case} $\bm{\mathcal{I}_t \subseteq \mathcal{I}_{c}}$: We rewrite $\widehat{\tensor{A}}_{\mathcal{I}_c}$ as an $\mathcal{I}_c$-mode multiplication using $\Pi_{\mathcal{I}_t}$ so that:
			\begin{align*}
				\sigma_k(\widehat{\tensor{A}}_{\mathcal{I}_c}) = \sigma_k((\Pi_{\mathcal{I}_t} \otimes I) P \tensor{A}_{\mathcal{I}_c}) \le \|(\Pi_{\mathcal{I}_t} \otimes I) P\|_2 \sigma_k(\tensor{A}_{\mathcal{I}_c}),
		\end{align*}
		where we used the property $\sigma_k(AB) \le \|A\|_2 \sigma_k(B)$~\cite{hackbusch_tensor_2019,mirsky_symmetric_1960} and $P$ permutes the factors of $\mathcal{I}_c$ so that $\mathcal{I}_t$ appears first. Then, using the orthogonal invariance of the spectral norm $\|.\|_2$, $\|A \otimes B\|_2 = \|A\|_2 \|B\|_2$, and $\|\Pi_{\mathcal{I}_t}\|_2 = 1 = \|I\|_2$, the proof is concluded.
		
		\textbf{Case }$\bm{\mathcal{I}_c \cap \mathcal{I}_{t} = \emptyset}$: We again rewrite $\widehat{\tensor{A}}_{\mathcal{I}_c}$ as a mode multiplication: 
				\[\sigma_k(\widehat{\tensor{A}}_{\mathcal{I}_c}) = \sigma_k(\tensor{A}_{\mathcal{I}_c} P (\Pi_{\mathcal{I}_t} \otimes I)^\top) = \sigma_k((\Pi_{\mathcal{I}_t} \otimes I) P^\top \tensor{A}_{\mathcal{I}_c}^\top ),\]
	where P is a permutation that permutes the factors of $\mathcal{I}^\perp_c := \mathcal{I} \backslash \mathcal{I}_c$ so that $\mathcal{I}_t \subseteq \mathcal{I}^\perp_c$ appears in the first positions. The proof is then completed as in the $\bm{\mathcal{I}_t \subseteq \mathcal{I}_{c}}$ case.
	\end{proof}
\end{lemma}
Using this interlacing property, a better bound on the error can be achieved. To apply M-HOSVD to the layerwise flattenings associated with a	dimension tree $\mathcal{T}$, we introduce notation for the corresponding layer partitions. For a layer $l$, let $\mathcal{T}_l = \{t_1,\dots,t_{m_l}\}$ and define the layer index set $S_l := \{1,\dots,m_l\}$. For each parent node $\ell \in \mathcal{T}_{l-1}$, we define its layer-$l$ children index block by $\mathcal{J}^{(l)}_\ell := \{ s \in S_l \mid t_s \in \operatorname{children}(\ell) \}$. The resulting layer-$l$ partition is then
	\[
	\mathfrak{J}_l := \bigl\{ \mathcal{J}^{(l)}_\ell \mid \ell \in \mathcal{T}_{l-1} \bigr\}.
	\]

\begin{theorem}
	\label{thm:RTL_error}
	Let $\tensor{A} \in \mathrm{V}_1 \otimes \dots \otimes \mathrm{V}_d$ be a tensor. The algorithm from \cref{def:RTL_projection} using M-HOSVD as the Tucker projection in each $P_l$ yields the error bound:
		\begin{equation*}
			\|\tensor{A} - (P_L \circ \cdots \circ P_1)(\tensor{A})\|
			\le \sum_{l=1}^{L} \|(P_l^\perp \circ P_{l-1} \circ \cdots \circ P_1) (\tensor{A})\|
			\le \|\tensor{A} - \tensor{A}^{\mathrm{best}}\| \sum_{l=2}^{L} \sqrt{|P_l|} ,
		\end{equation*}
		where $|P_l|$ denotes the number of SVD projections on layer $l$, $\tensor{A}^{\mathrm{best}}$ is the best HTD approximation with the same hierarchical rank at each node and $P_l^\perp := I-P_l$.
	\begin{proof}
		We show that $\|(P_l^\perp ((P_{l-1} \circ \cdots \circ P_1) (\tensor{A}))\| \le \|P_l^\perp (\tensor{A})\|$. Using an M-HOSVD projection as the Tucker projection in $P_l$ with grouping $\mathfrak{J}_l = \{ \mathcal{J}_{\ell_1}, \dots, \mathcal{J}_{\ell_n}\}$ corresponds to $n$ ST-HOSVD projections $P_{\mathcal{J}_{\ell_1}} \boxtimes \dots \boxtimes P_{\mathcal{J}_{\ell_n}}$, where each $P_{\mathcal{J}_\ell}$ is applied to $\tensor{X}_{\mathcal{T}_l} := ((P_{l-1} \circ \dots \circ P_{1})\tensor{A})_{\mathcal{T}_l}$. In total $m$ SVD projections $P_{\mathcal{I}_{t_i}}$ are performed, where $\mathcal{T}_l = \{t_1,\dots,t_m\}$. Then by the proof of \cref{lemma:MTHOSVD_err},
		\[
			\|P_l^\perp \tensor{X}\|^2 \le \sum_{i =1}^n \|P_{\mathcal{J}_{\ell_i}}^\perp \tensor{X}\|^2 \le \sum_{t \in \mathcal{T}_l} \|P_{\mathcal{I}_t}^\perp \tensor{X}\|^2.
		\]
		In a dimension tree, for any two nodes $t$ and $k$, either $\mathcal{I}_t \subseteq \mathcal{I}_k$, $\mathcal{I}_k \subseteq \mathcal{I}_t$ or $\mathcal{I}_t \cap \mathcal{I}_k = \emptyset$. Each previous layerwise projection $P_i$ for $i <l$, can be written as a sequence of SVD projections $P_{\mathcal{I}_k}$ for $k \in \mathcal{T}_i$. So, due to \cref{lemma:interlacing_singval}, $\sigma(\tensor{X}_{\mathcal{I}_t}) \le \sigma(\tensor{A}_{\mathcal{I}_t})$, componentwise, and as a result $\|P_{\mathcal{I}_t}^\perp \tensor{X}\| \le \|P_{\mathcal{I}_t}^\perp \tensor{A}\|$, concluding the proof.
	\end{proof}
\end{theorem}
\begin{remark}
	\label{remark:error_orthogonality}
	\cref{thm:RTL_error} can be sharpened further when the terms in the telescoping sum are orthogonal. This occurs when the projections commute, which happens when each $P_l$ projects orthogonally onto a subspace of $\operatorname{Im}(P_k)$ for all layers $k <l $. Projectors computed in the LTR direction satisfy this last property, explaining why the LTR algorithm admits the sharper error bound as shown in \cite{hackbusch_tensor_2019}:
	\[
	\|\tensor{A} - \tensor{A}^{\mathrm{LTR}}\|
	\le \sqrt{2d - 3} \cdot \|\tensor{A} - \tensor{A}^{\mathrm{best}}\|.
	\]
\end{remark}

\begin{remark}
	Trees with the same amount of nodes, but different heights or slicings into layers, lead to different error bounds for RTL. Since $L \le d$ for any dimension tree, we can, however, give an upper bound on the factor $\sum_{l=2}^{L} \sqrt{|P_l|}$. Using the Cauchy--Schwarz inequality $\sum_{l=2}^{L} \sqrt{|P_l|} \le \sqrt{(L-1)\sum_{l=2}^{L} |P_l|}$, $L \le d$, and $\sum_{l=2}^{L} |P_l| \le 2d-2$ for any dimension tree, we find $\sum_{i=2}^L \sqrt{|P_i|} \le \sqrt{2}(d-1)$.
	This leads to an error upper bound over all possible tree structures equal to
	\[\|\tensor{A} - (P_L \circ \cdots \circ P_1)(\tensor{A})\| \le \sqrt{2}(d-1)\|\tensor{A} - \tensor{A}^{\mathrm{best}}\|,\]
	showing that the error bound of RTL with M-HOSVD projections will at most be a factor $\sqrt{d}$ (for $d \ge 2$) larger than the error bound of LTR.
\end{remark}

\subsection{An efficient truncation algorithm into SVD-HTD form}
\label{sec:algorithm}
A naive realization of the projection scheme described in \cref{def:RTL_projection} consists of reconstructing the full tensor $\tensor{A}$ after each layer and performing an M-HOSVD projection on $\tensor{A}_{\mathcal{T}_l}$ at the next layer $l$. This approach is computationally prohibitive for large-scale tensors.

Next, we exploit the SVD-HTD form to substantially reduce the complexity of the M-HOSVD projections, by using the following lemma and its corollary.
\begin{lemma}
	\label{lemma:approx_Q}
	Let $\tensor{A} \in \mathrm{V}_1 \otimes \dots \otimes \mathrm{V}_d \otimes \mathrm{E}$ be a tensor. Let $P$ denote an M-HOSVD approximation of the form $P = P_{\mathcal{J}_1} \boxtimes \dots \boxtimes P_{\mathcal{J}_n} \otimes I$, using a partition $\mathcal{J}_1 \sqcup \dots \sqcup \mathcal{J}_n = (1,\dots, d)$, i.e., no truncation in the factor $\mathrm{E}$. If $Q \in \mathbb{R}^{N \times \dim(\mathrm{E})}$ is a matrix with orthonormal columns, then, $Q \cdot_{d+1} P(\tensor{A}) = P(Q \cdot_{d+1} \tensor{A})$.
		\begin{proof}
		Define $\widetilde{\tensor{A}} := Q \cdot_{d+1} \tensor{A}$. For any mode-$k$ flattening with $k \le d$, we have
		\begin{align*}
			\widetilde{\tensor{A}}_{(k)} &= \tensor{A}_{(k)}(I \otimes Q)^\top,
		\end{align*}
		where $I$ is a suitably sized identity matrix. Since $(I \otimes Q)^\top$ has orthonormal rows, $\tensor{A}_{(k)}$ and $\widetilde{\tensor{A}}_{(k)}$ share left singular subspaces and singular values, concluding the proof.
	\end{proof}
\end{lemma}

\begin{corollary}
	\label{cor:HOSVD_child}
	Let the SVD-HTD property hold at layer $l$ for the HTD $(\mathcal{T},\bm{B})$ representing 
	$\tensor{A}$. Let $P_{l+1} = P_{\mathcal{J}_{\ell_1}} \boxtimes \dots \boxtimes P_{\mathcal{J}_{\ell_n}}$ be an M-HOSVD projection with $\mathfrak{J}_l = \{\mathcal{J}_{\ell_1},\dots,\mathcal{J}_{\ell_n}\}$. For each node $\ell \in \mathcal{T}_l$ with $\operatorname{children}(\ell) = \{t_1, \dots, t_m\}$, let $F_{\ell} $ denote the frame defined recursively by \cref{eq:rec_def}. Then $P_{l+1}(\tensor{A}_{\mathcal{T}_{l+1}})$ can be computed by evaluating the ST-HOSVD projections $P_{\{1,\dots,m\}}((F_\ell)_{\{\mathcal{I}_{t_1},\dots,\mathcal{I}_{t_m}, e_\ell\}})$ for all $\ell \in \mathcal{T}_{l}$.
	\begin{proof}
		The Tucker representation of $\tensor{A}_{\mathcal{T}_l}$ equals $(F_{\ell_1},\dots, F_{\ell_n})\cdot\tensor{G}_l$, where $F_{\ell_k}$ is derived from \cref{eq:rec_def}. The $\mathcal{I}_{\ell_k}$-flattening of $\tensor{A}$ can be expressed as
		\[
		\tensor{A}_{\mathcal{I}_{\ell_k}} = F_{\ell_k} (\tensor{G}_l)_{(k)}(F_{\ell_1} \otimes \dots \otimes F_{\ell_{k-1}} \otimes F_{\ell_{k+1}} \otimes \dots \otimes F_{\ell_n})^\top.
		\] By \cref{def:svd_htd}, the SVD-HTD property ensures that $F_{\ell_k} = U_{\ell_k} \Sigma_{\ell_k}$, where $\tensor{A}_{{\ell_k}} = U_{\ell_k} \Sigma_{\ell_k} V_{\ell_k}^\top$ is an SVD. Due to the equality of both expressions, we conclude that 
		\[V_{\ell_k}^\top = (\tensor{G}_l)_{(k)}(F_{\ell_1} \otimes \dots \otimes F_{\ell_{k-1}} \otimes F_{\ell_{k+1}} \otimes \dots \otimes F_{\ell_n})^\top. \]
		The statement then follows from applying \cref{lemma:approx_Q} with $Q = V_{\ell_k}$ for each $\ell_k$.
	\end{proof}
\end{corollary}

We now present a practical and efficient implementation of the RTL algorithm that leverages the sequential error bound established in \cref{thm:RTL_error}. We propose \cref{alg:RTL_SVD_HTD} for approximating a general dense tensor by an SVD-HTD using an RTL approach. Crucial for the computational efficiency of this scheme is the use of randomized SVDs in line~5. The main loop of the algorithm is conceptually analogous to the second stage of \cref{alg:SVD_HTD}, with the exception that truncating ST-HOSVDs are used in place of their non-truncating counterparts. To conceptualize this equivalence, one may assume an HTD of $\tensor{A}$ exists for which every non-root $B_t$ is a suitably sized identity matrix, while $B_1 = \operatorname{vec}(\tensor{A})$.

The correctness of \cref{alg:RTL_SVD_HTD} is established in the next proposition.
\begin{algorithm}[tb]
	\caption{RTL SVD-HTD}
	\begin{algorithmic}[1]
	\label{alg:RTL_SVD_HTD}
		\REQUIRE Tensor $\tensor{A} \in \mathrm{V}_1 \otimes \dots \otimes \mathrm{V}_d$, dimension tree $\mathcal{T}$, hierarchical ranks $r_t$ for each node $t \in \mathcal{T}$.\\
		\STATE $B_1 \gets \operatorname{vec}(\tensor{A})$\\
		\FOR{$l \gets 1$ to $L-1$}
			\FOR{$t \in \mathcal{T}_l \backslash (\bigcup_{j<l} \mathcal{T}_j \cup \mathcal{T}_L)$}
				\STATE $\{t_1, \dots ,t_m\} \gets \operatorname{children}(t)$\\
				\STATE Let $(\tensor{B}_t)_{\{\mathcal{I}_{t_1},\dots, \mathcal{I}_{t_m}, e\}} \approx (\widehat{U}_{t_1},\dots, \widehat{U}_{t_m}, U_e)\cdot \tensor{G}$ be the ST-HOSVD of $B_t$, flattened to the index sets of $t$'s children. That is, \\$(B_t)_{\mathcal{I}_{t_k}} = U_{t_k} \Sigma_{t_k} V_{t_k}^\top = \begin{bmatrix}\widehat{U}_{t_k} & \widetilde{U}_{t_k}\end{bmatrix}
				\begin{bmatrix}\widehat{\Sigma}_{t_k} & \\ & \widetilde{\Sigma}_{t_k}\end{bmatrix}
				\begin{bmatrix}\widehat{V}^\top_{t_k} \\ \widetilde{V}^\top_{t_k}\end{bmatrix}$.\\
				\STATE $(\tensor{B}_t)_{\{\mathcal{I}_{t_1},\dots, \mathcal{I}_{t_m}, e\}} \gets (\widehat{\Sigma}_{t_1}^{-1} ,\dots, \widehat{\Sigma}_{t_m}^{-1}, U_e) \cdot \tensor{G}$\\
				\STATE $B_{t_j} \gets \widehat{U}_{t_j} \widehat{\Sigma}_{t_j}$ 		for $j = 1,\dots, m$\\
			\ENDFOR\\
			\COMMENT{SVD-HTD Normalization}\\
			\STATE Run \cref{alg:SVD_HTD} on the subtree of $\mathcal{T}$ of depth $l+1$, rooted at the same root as $\mathcal{T}$, using the associated $B_t$.\\
		\ENDFOR\\
		\RETURN $(\mathcal{T}, \bm{B})$
	\end{algorithmic}
\end{algorithm}

\begin{proposition}
	\cref{alg:RTL_SVD_HTD} computes an RTL approximation of the input tensor $\tensor{A} \in \mathrm{V}_1 \otimes \dots \otimes \mathrm{V}_d$, as in \cref{def:RTL_projection} where each $P_l$ uses an M-HOSVD projection. The output is in SVD-HTD form. 
	\begin{proof}
		At the beginning of each new layer $l$ in line~2, the HTD $(\mathcal{T},\bm{B})$ of the approximation $(P_{l-1} \circ \dots \circ P_1)(\tensor{A})$ is in SVD-HTD form for all layers $1, \dots,l$. This holds either by the initialization $B_1 = \operatorname{vec}(\tensor{A})$ in line~1 or due to line~9 of the previous iteration. Consequently, \cref{cor:HOSVD_child} implies that the ST-HOSVDs computed in line~5 on all nodes in layer $l$ are equivalent to performing an M-HOSVD on the flattening of the represented tensor at that point $((P_{l-1} \circ \dots \circ P_1)(\tensor{A}))_{\mathcal{T}_l}$ using the partitioning $\mathfrak{J}_l$. Note that each $B_t$ at the beginning of a loop in line 3, is the $F_t$ for that node, as is evident by its assignment in line~7, and only becomes the transfer tensor after being overwritten in line~6. Each iteration of line~2 therefore iteratively builds the sequential projection chain from \cref{def:RTL_projection}. After each M-HOSVD, however, the SVD-HTD form is lost due to truncation and must be restored in line~8 before proceeding to the next layer.
\end{proof}
\end{proposition}

This result implies in particular that the error bound from \cref{thm:RTL_error} applies. Hence, \cref{alg:RTL_SVD_HTD} yields a quasi-optimal approximation with a quasi-optimality factor $\sqrt{2}d$. This is slightly worse than the factor $\sqrt{2d - 3}$ of the LTR approximation from \cref{def:LTR_projection}, but in practice, the factor of RTL is observed to closely adhere to the $\sqrt{2d -3}$ in most cases; see \cref{sec:experiments}.

\begin{remark}
	\label{rem:simple_RTL}
	Although line 9 of \cref{alg:RTL_SVD_HTD} does not change the overall complexity of the algorithm, as it only involves the reduced transfer tensors $\tensor{B}_t$; nevertheless, for efficiency, the SVD-HTD normalization step can be omitted. 
	The correctness of this variant, follows from~\cite[Theorem 11.63]{hackbusch_tensor_2019}. While this theorem applies to a procedure for truncating tensors already in HTD form, we observe that setting $B_t = \tensor{A}_{(\mathcal{I}_l, \mathcal{I}_r)}$ at the root node $t$ with children $l$ and $r$, and $B_{t'} = I$ at all other nodes of the dimension tree, extends its applicability to general dense tensors.
	Overall, we think the error analysis and proof of correctness are simpler with our approach.
	At the $l$th layer, the algorithms with SVD-HTD form, compute the singular values of $(P_{l-1} \circ \dots \circ P_1)(\tensor{A})$ instead of the overestimation computed when SVD-HTD form is not used. Therefore, the former should more accurately estimate the approximation error, which in turn can increase the compression factor, as we will see in \cref{sec:experiments}.
\end{remark}

 \begin{remark}	
 	To control the truncation accuracy during the RTL procedure, we use the middle bound of \cref{thm:RTL_error}, which decomposes the error into layerwise contributions $\|(P_l^\perp \circ P_{l-1} \circ \cdots \circ P_1)(\tensor{A})\|$. These quantities determine how much error is introduced by the truncation at layer $l$. A common strategy is to distribute a total tolerance $\varepsilon$ evenly across all layers, e.g., by assigning each layer an error budget $\varepsilon_l = \varepsilon / L$. In the sequential setting, we further allow unused error from earlier nodes, i.e., the difference between the	allocated and the estimated error, to be passed forward and added to the budgets of subsequent nodes, ensuring that the full error budget is utilized. By the quasi-optimality of the M-HOSVD, we have
 		\[
 		\|(P_l^\perp \circ P_{l-1} \circ \cdots \circ P_1)(\tensor{A})\|^2 
 		\;\le\;
 		\sum_{t_j \in \mathcal{T}_l}
 		\sum_{i > r_{t_j}} \sigma_{i,t_j}^2 ,
 		\]
 		where $\sigma_{i,t_j}$ is the $i$th singular value of $((P_{l-1} \circ \cdots \circ P_1)(\tensor{A}))_{\mathcal{T}_{t_j}}$ and $r_{t_j}$ is the truncation rank for every node $t_j \in \mathcal{T}_l$. Hence, truncation at layer $l$ may be performed by selecting the truncation rank $r_{t_j}$ as the smallest integer such that
 		\[
 		\sum_{t_j \in \mathcal{T}_l}
 		\sum_{i > r_{t_j}} \sigma_{i,t_j}^2 \;\le\; \varepsilon_l^2 .
 		\]
 		In principle, the interaction of errors between layers should be included to match the theoretical upper bound exactly. However, simply summing the squared discarded singular values at each layer, analogous to the standard LTR procedure, always accurately bounded the approximation error in our experiments.
 \end{remark}
\subsection{Complexity}
\label{sec:complexity}

The asymptotic computational complexity of~\cite[Algorithm 2]{grasedyck_hierarchical_2010} and \cref{alg:RTL_SVD_HTD} can be bounded in terms of the maximal dimension size $N = \max \{n_1, \dots, n_d\}$ and the maximal rank at the $l$th layer $R_{(l)} = \max_{t\in \mathcal{T}_l} r_t$. The dominant cost originates from the truncated SVDs. We assume the asymptotic time complexity for computing a rank-$r$ truncated randomized SVD of an $m\times n$ matrix is $\mathcal{O}(mnr)$. Then, after some straightforward counting, the complexities for computing the HTD of a dense tensor using a perfect binary tree can be found in \cref{tab:complexity_balanced}.

A key limitation of the LTR direction is the appearance of large exponents of $R_{(l)}$ in the complexity expressions. This is particularly detrimental for HTD, as the hierarchical rank $R_{(l)}$ is not necessarily bounded by $N$. The RTL direction, on the other hand, suffers from the fact that the hierarchical rank can increase (possibly exponentially) towards the root. This means that the first SVD, which is the most expensive, could be much more expensive than the first SVD of the LTR direction. Naturally, if $R_{(2)} \le R_{(L)}$, RTL should outperform LTR, as the exponents of $N$ decay quickly after the first SVD. Nevertheless, as shown in \cref{sec:experiments}, even in cases where $R_{(2)} > R_{(L)}$, RTL can still yield significant speedups in practice.

\begin{table}[tb]
	\centering
	\footnotesize
	\caption{\footnotesize Asymptotic computational complexity of the LTR algorithm of Grasedyck~\cite[Algorithm 2]{grasedyck_hierarchical_2010} and our RTL algorithm, \cref{alg:RTL_SVD_HTD}, for a perfect binary tree structure (so $L = \log_{2}(d)$) using T-HOSVD and ST-HOSVD.}\label{tab:complexity_balanced}
	
	\begin{tabular}{llc}
		\toprule
		\textbf{Tucker} & \textbf{Direction} 
		& \textbf{Asymptotic complexity} \\
		\midrule
		\multirow{2}{*}{T-HOSVD} 
		& LTR & $d N^d R_{(L)} + \sum_{l=3}^{L-1} \left(2^{l-1} R^{2^{l}}_{(l+1)} R_{(l)}\right) + R_{(3)}^4 R_{(2)}$ \\[4pt]
		& RTL & $N^d R_{(2)} + \sum_{l = 2}^{L-1} 2^{l} N^{\frac{d}{2^{l-1}}}R_{(l)}R_{(l+1)}$ \\[4pt]
		\midrule
		\multirow{2}{*}{ST-HOSVD}
		& LTR & $\sum_{i=1}^{d} N^{d-i+1} R_{(L)}^{i} + \sum_{l=3}^{L-1}\sum_{i=1}^{2^{l-1}} \left(R^{2i}_{(l+1)} R^{2^{l-1}-i+1}_{(l)} \right)+ R^4_{(3)} R_{(2)}$ \\[4pt]
		& RTL & $N^d R_{(2)} + \sum_{l = 2}^{L-1} 2^{l-1} \left(N^{2^{i}} R_{(l)}R_{(l+1)} + N^{2^{i-1}} R_{(l)}R_{(l+1)}^2\right)$ \\[4pt]
		\bottomrule
	\end{tabular}
\end{table}

The dimension tree heavily influences the computational complexity of the RTL algorithm. Trees which do not generate approximately square matricizations of the tensor, inhibit the exponential decay of the exponent of $N$ for the RTL direction, although the ST-HOSVD method reduces the effects of this interference.

\begin{figure}[tb]
	\centering
	\begin{tikzpicture}
		\begin{axis}[
			width=0.9\textwidth,
			height=0.4\textwidth,
			xmode=log,
			ymode=log,
			grid=both,
			grid style={dotted},
			xlabel={Aspect ratio $m/n$},
			ylabel={Throughput [GFLOP/s]},
			xmin=1e-6,
			xmax=1e6,
			ymin=1,
			ymax=300,
			tick align=outside,
			xtick={1e-6,1e-4,1e-2,1,1e2,1e4,1e6},
			legend style={
				draw=none,
				fill opacity=0.6,
				text opacity=1,
				font=\footnotesize,
				at={(0.03,0.97)},
				anchor=north west
			},
			legend cell align=left
			]
			
			\addplot[
			thick,
			color=blue!70!black,
			mark=square*,
			mark options={fill=blue!70!black}
			] coordinates {
				(1e-6, 2.943) (1e-4, 7.359) (1e-2, 85.76) (1, 61.30) (1e2, 122.7) (1e4, 12.81) (1e6, 4.525)
			};
			\addlegendentry{Deterministic}
			
			\addplot[
			thick,
			color=orange!85!black,
			mark=*,
			mark options={fill=orange!85!black}
			] coordinates {
				(1e-6, 4.151) (1e-4, 16.27) (1e-2, 62.00) (1, 92.92) (1e2, 55.08) (1e4, 18.40) (1e6, 4.215)
			};
			\addlegendentry{Randomized}
			
			\draw[dashed, thick] (axis cs:1,1) -- (axis cs:1,300);
			
		\end{axis}
	\end{tikzpicture}
	\caption{\footnotesize
		Throughput versus matrix aspect ratio for the rank-$r$ SVD approximation of an $m\times n$ matrix with $10^8$ elements, where $r = \min(m,n,500)$.
	}
\label{fig:svd_throughput}
\end{figure}
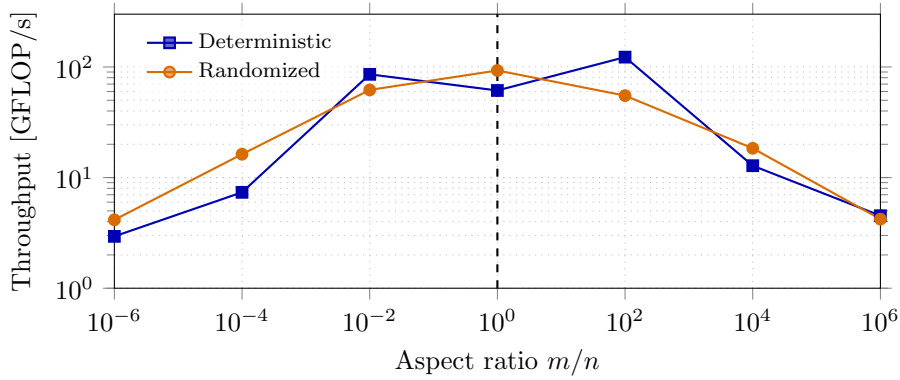

 Finally, the throughput of different matrix shapes also influences timing results. \Cref{fig:svd_throughput} illustrates that for both deterministic and randomized SVDs, throughput peaks near square matrices
 and degrades symmetrically for thin and wide matrices. This further explains why RTL can outperform LTR, even when the hierarchical ranks near the root are larger. For instance, for the timings of \cref{fig:svd_throughput}, computing a rank-$500$ approximation of a $10\,000\times10\,000$ matrix using a randomized SVD was $1.5 \times$ faster than computing a rank-$10$ approximation of a $10\times 10\,000\,000$ matrix even though the asymptotic complexity analysis suggests that it should be $50$ times slower. 

\section{Numerical experiments}
\label{sec:experiments}
The different methods described above are tested on three different test cases to demonstrate their applicability and performance.
To demonstrate the scaling properties of the proposed methods, they are applied to different sizes of the Hilbert tensor, a toy benchmark problem.
We also evaluate the proposed HTD schemes on two realistic, large-scale, high-dimensional datasets. 
Each dataset was considered under two tensorizations:
\begin{enumerate}
	\item A \textbf{moderate-order} form corresponding to the natural structure of the data.
	\item An \textbf{order-16} form obtained through tensorization, i.e., splitting modes into smaller factors, increasing the order of the tensor while preserving the total number of entries, e.g., $\mathbb{R}^{10} \simeq \mathbb{R}^{5\times 2}$.
\end{enumerate}
This allows us to examine how performance is influenced by the tensor order on practical examples.

There are three methods that we will test: \begin{itemize}
	\item LTR: The default LTR algorithm as defined by Grasedyck in~\cite[Algorithm 2]{grasedyck_hierarchical_2010} and described in \cref{sec:htd-comp}. This is the algorithm used in libraries that compute the HTD of a tensor~\cite{perisa_tensortoolboxjl_2021,kressner_algorithm_2014}.
	\item RTL: \cref{alg:RTL_SVD_HTD}, but without the SVD-HTD transformation in line 8. This corresponds to a general dense tensor version of the sequential RTL algorithm of Hackbusch~\cite{hackbusch_tensor_2019} but using randomized truncated SVDs.
	\item M-HOSVD RTL: \cref{alg:RTL_SVD_HTD}.
\end{itemize}
Each method will be tested using either sequential or parallel truncation inside each layer. That is, for LTR the T-HOSVD in \cref{def:LTR_projection} is replaced by an ST-HOSVD in the sequential case. For the RTL methods, we use ST-HOSVD and the T-HOSVD, respectively, for the Tucker approximation in lines 5--6 of \cref{alg:RTL_SVD_HTD}. For RTL the overall approximation in each layer is an M-HOSVD when ST-HOSVD is used.

Crucial for the performance of all methods, in particular those using the RTL direction, is that randomized truncated SVD are used, rather than truncating a compact SVD. Each method uses the same one-sided sketch randomized SVD with either a column or row space range finder \cite{halko_finding_2011}. The decision to use a column or row space range finder is decided per matricization to minimize the size of the matrix on which the deterministic SVD is performed.
Randomized SVDs using sketching must define the dimension of the sketching matrix, $\bm{\Omega}$. For this rank, an initial guess $R^*$ is chosen for all nodes and an oversampling rate of $5$, so that the initial $\bm{\Omega} \in \mathbb{R}^{n \times (R^* + 5)}$. A blocked version of the adaptive randomized range finder~\cite[Algorithm 4.2]{halko_finding_2011} is then used to iteratively grow the orthonormal basis by adding $R^*$ columns each iteration. This expansion is continued until the norm of the residual after orthogonalization, against the previous basis $Q$, of a new $\bm{A}\bm{\Omega}$ is small enough in norm, i.e., $\|(I- QQ^\top) \bm{A}\bm{\Omega}\| \le 10 \sqrt{\frac{2}{\pi}} \epsilon$, where $\epsilon$ is the desired error threshold. Using block increments of $R^*$ is more computationally efficient than increasing the basis by one vector each iteration, but can lead to an overestimation of the rank. We do, however, bound the used sketching dimension by the minimal dimension of the matrix.

All tensors were normalized ($\|\tensor{A}\| = 1$) and then compressed with a desired relative Frobenius error of $\varepsilon = 10^{-2}$ and a balanced binary dimension tree. For all methods, the squared error is estimated as the sum of the squared discarded singular values at each node, this guides the truncation ranks ensuring the total error is less than $\varepsilon$. Note this assumes the errors between layers are orthogonal to each other, which is in contrast to \cref{thm:RTL_error} for RTL. Nonetheless, even with this inaccuracy, we will see that none of the RTL errors exceed the desired error bound. We do not explicitly report the estimated error, since for all experiments the estimated error bounded the computed error, while agreeing in the first two significant digits.
The error of the approximated tensor $\widehat{\tensor{A}}$ against the original $\tensor{A}$ is measured using: $\frac{\|\tensor{A} - \widehat{\tensor{A}}\|_F}{\|\tensor{A}\|_F}$. The compression factor is computed as $\frac{|\tensor{A}|}{\sum_{t \in \mathcal{T}} |B_t|}$, where $|.|$ denotes the number of elements in the coordinate representation of the argument.

We implemented these method in Julia version 1.11.5~\cite{bezanson_julia_2017}, although we also tested our implementation on Julia version 1.12.1. Our code depends on the LinearAlgebra and Random libraries. The source code needed to run the experiments is publicly available\footnote{\url{https://gitlab.kuleuven.be/numa/public/m-hosvd-htd}}. We measure the execution time and total memory allocation of our implementations using the BenchmarkTools.jl library~\cite{revels_benchmarktoolsjl_2025}. The underlying LAPACK version is 3.12.0 which calls the OpenBLAS library using the ILP64 index bindings. No explicit parallelization was used besides the inherent parallelism of the LAPACK and BLAS calls.
The Julia experiments were run on the wICE cluster of the Flemish Supercomputer Centre (VSC) on a single IceLake node (Intel Xeon Platinum 8360Y CPUs@2.4 GHz with 54MB of L3 cache) using 36 cores. 

\subsection{Hilbert tensor}
\label{sec:hilb}
To validate the scaling properties predicted by our complexity analysis empirically, we first test the different algorithms on the Hilbert tensor
\[
\tensor{A}_{i_1, \dots, i_d} = \frac{1}{1 - d + \sum_{k=1}^d i_k}, \quad 1 \le i_k \le N,\; \; k = 1,\dots,d,
\]
using varying orders $d$ and sizes $N$. Due to the rapid decay of singular values, the hierarchical ranks required for an accurate HTD approximation are tiny, typically $ \max_{t\in \mathcal{T}} r_t \approx 3$ in our experiments. The initial sketching rank was set to $R^* = 20$.

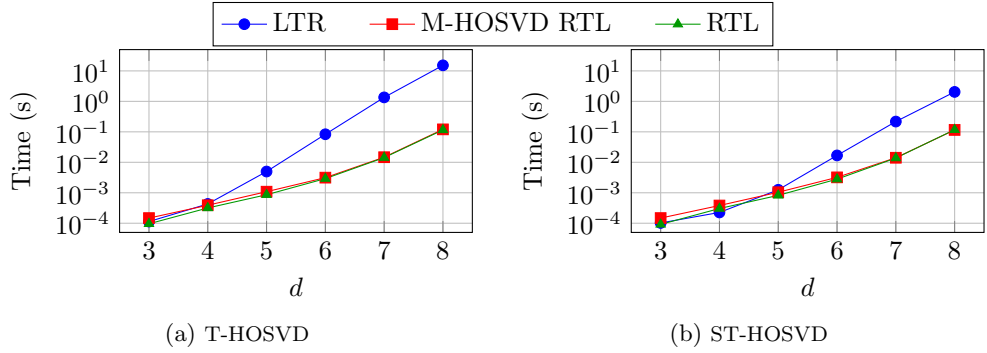
\begin{figure}[tb!]
	\centering
	\pgfplotstableread[col sep=comma]{./anc/timings_N_10.csv}\datatable
	\pgfplotscreateplotcyclelist{colors}{
		{blue,mark=*},
		{red,mark=square*},
		{green!60!black,mark=triangle*},
		{black,mark=o},
	}
\ref{sharedlegend}\\
	\begin{subfigure}{0.48\textwidth}
		\centering
		\begin{tikzpicture}
			\begin{axis}[
				xmode=normal,
				ymode=log,
				xlabel={$d$},
				ylabel={Time (s)},
				width=\linewidth,
				height=4cm,
				ytick style={draw=none},
				grid=major,
				xtick=data,
				ytick = {1e-4,1e-3,1e-2,1e-1,1e0,1e1,1e2},
				ymin = 5e-5,
				ymax = 5e1,
				mark options=solid,
				cycle list name=colors,
				legend entries={LTR, M-HOSVD RTL, RTL},
				legend to name=sharedlegend,
				legend columns=-1, 
				legend style={
					/tikz/every even column/.append style={column sep=0.5cm}
				},
				]
				\addplot table[x=d,y=HTUCKER_LTR]{\datatable};
				
				\addplot table[x=d,y=adapted_HTUCKER_RTL]{\datatable};
				
				\addplot table[x=d,y=HTUCKER_RTL]{\datatable};
			\end{axis}
		\end{tikzpicture}
		\caption{\footnotesize T-HOSVD}
	\end{subfigure}
	\hfill
	\begin{subfigure}{0.48\textwidth}
		\centering
		\begin{tikzpicture}
			\begin{axis}[
				xmode=normal,
				ymode=log,
				xlabel={$d$},
				ylabel={Time (s)},
				ytick style={draw=none},				
				ytick = {1e-4,1e-3,1e-2,1e-1,1e0,1e1,1e2},
				ymin = 5e-5,
				ymax = 5e1,
				width=\linewidth,
				height=4cm,
				grid=major,
				xtick=data,
				mark options=solid,
				cycle list name=colors
				]
				\addplot table[x=d,y=ST_HTUCKER_LTR]{\datatable};
				
				\addplot table[x=d,y=adapted_ST_HTUCKER_RTL]{\datatable};

				\addplot table[x=d,y=ST_HTUCKER_RTL]{\datatable};
			\end{axis}
		\end{tikzpicture}
		\caption{\footnotesize ST-HOSVD}
	\end{subfigure}
	\caption{\footnotesize Timing comparison of LTR vs.\ RTL for the Hilbert tensor for different orders $d$, where each dimension has size $N = 10$.}
	\label{fig:d-timings}
\end{figure}

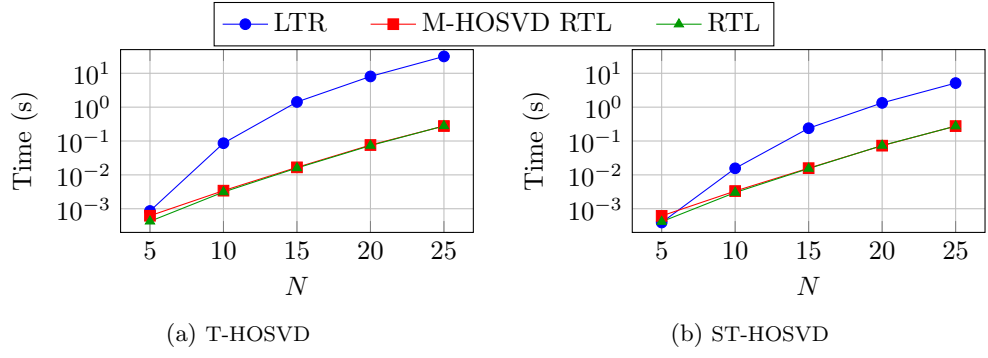
\begin{figure}[tb!]
	\centering
	\pgfplotstableread[col sep=comma]{./anc/timings_d_6.csv}\datatable
	\pgfplotscreateplotcyclelist{colors}{
		{blue,mark=*},
		{red,mark=square*},
		{green!60!black,mark=triangle*},
		{black,mark=o},
	}\ref{sharedlegend}\\

	\begin{subfigure}{0.48\textwidth}
		\centering
		\begin{tikzpicture}
			\begin{axis}[
				xmode=normal,
				ymode=log,
				xlabel={$N$},
				ylabel={Time (s)},
				width=\linewidth,
				height=4cm,
				ytick style={draw=none},
				grid=major,
				xtick=data, 
				ytick = {1e-3,1e-2,1e-1,1e0,1e1},
				ymin = 2e-4,
				ymax = 5e1,
				mark options=solid,
				cycle list name=colors,
				legend entries={LTR,  M-HOSVD RTL, RTL},
				legend to name=sharedlegend,
				legend columns=-1, 
				legend style={
					/tikz/every even column/.append style={column sep=0.5cm}
				},
				]
				\addplot table[x=N,y=HTUCKER_LTR]{\datatable};
				
				\addplot table[x=N,y=adapted_HTUCKER_RTL]{\datatable};
				
				\addplot table[x=N,y=HTUCKER_RTL]{\datatable};
			\end{axis}
		\end{tikzpicture}
		\caption{\footnotesize T-HOSVD}
	\end{subfigure}
	\hfill
	\begin{subfigure}{0.48\textwidth}
		\centering
		\begin{tikzpicture}
			\begin{axis}[
				xmode=normal,
				ymode=log,
				xlabel={$N$},
				ylabel={Time (s)},
				ytick style={draw=none},
				width=\linewidth,
				height=4cm,
				grid=major,
				xtick=data,				
				ytick = {1e-3,1e-2,1e-1,1e0,1e1},
				ymin = 2e-4,
				ymax = 5e1,
				mark options=solid,
				cycle list name=colors
				]
				\addplot table[x=N,y=ST_HTUCKER_LTR]{\datatable};
				
				\addplot table[x=N,y=adapted_ST_HTUCKER_RTL]{\datatable};

				\addplot table[x=N,y=ST_HTUCKER_RTL]{\datatable};
			\end{axis}
		\end{tikzpicture}
		\caption{\footnotesize ST-HOSVD}
	\end{subfigure}
	\caption{\footnotesize Timing comparison of LTR vs.\ RTL for the 6th-order Hilbert tensor for different dimension sizes $N$.}
	\label{fig:N-timings}
\end{figure}

The results of varying $d$ are shown in \cref{fig:d-timings}, while \cref{fig:N-timings} show the timings of varying $N$.
As anticipated, \cref{fig:d-timings,fig:N-timings} demonstrate that the RTL algorithms exhibit nearly identical runtime, reflecting their equal asymptotic time complexity for both T-HOSVD and ST-HOSVD variants. In contrast, there is a pronounced difference between the LTR and RTL approaches: for the ST-HOSVD case, the RTL algorithm is approximately one order of magnitude faster, while for the T-HOSVD variant, the speedup factor increases to about two orders of magnitude for larger values of $d$ and $N$. This is to be expected because $R_{(2)} = R_{(L)}$. Note that \cite[Algorithm 2]{grasedyck_hierarchical_2010}, found in current software libraries, corresponds to the T-HOSVD LTR algorithm.

\begin{figure}[tb!]
	
	\pgfplotstableread[col sep=comma]{./anc/STTNN_vs_STRTL.csv}\datatable
	\pgfplotscreateplotcyclelist{colors}{
		{blue,mark=*},
		{red,mark=square*},
		{green!60!black,mark=triangle*},
		{black,mark=o},
	}
	\centering
	\ref{sttnnlegend}\\
	\vfill
	\begin{subfigure}{0.48\textwidth}
		\centering
		\begin{tikzpicture}
			\begin{axis}[
				xmode=normal,
				ymode=log,
				xlabel={Rank},
				ylabel={Time (s)},
				width=0.9\linewidth,
				height=4cm,
				ytick style={draw=none},
				grid=major, 
				mark options=solid,
				cycle list name=colors,
				legend entries={STTNN, ST-HOSVD RTL},
				legend to name=sttnnlegend,
				legend columns=-1, 
				legend style={
					/tikz/every even column/.append style={column sep=0.5cm}
				},
				]
				\addplot table[x=rank,y=STTN]{\datatable};
				
				\addplot table[x=rank,y=STRTL]{\datatable};
			\end{axis}
		\end{tikzpicture}
		\caption{\footnotesize Timing}
	\end{subfigure}
	\hfill
	\begin{subfigure}{0.48\textwidth}
		\centering
		\begin{tikzpicture}
			\begin{axis}[
				xmode=normal,
				ymode=log,
				xlabel={Rank},
				ylabel={$\frac{\|\tensor{A} - \widehat{\tensor{A}}\|}{\|\tensor{A}\|}$},
				ytick style={draw=none},
				width=0.9\linewidth,
				height=4cm,
				grid=major,
				mark options=solid,
				cycle list name=colors
				]
				\addplot table[x=rank,y=STTN_Error]{\datatable};
				
				\addplot table[x=rank,y=STRTL_Error]{\datatable};
			\end{axis}
		\end{tikzpicture}
		\caption{\footnotesize Error}
	\end{subfigure}
	\caption{\footnotesize Timing and error comparison of RTL with ST-HOSVD and STTNN from \cite{bucci_randomized_2024} for the Hilbert tensor with $d =6$ and $N=20$.}
	\label{fig:sttnn-strtl}
\end{figure}
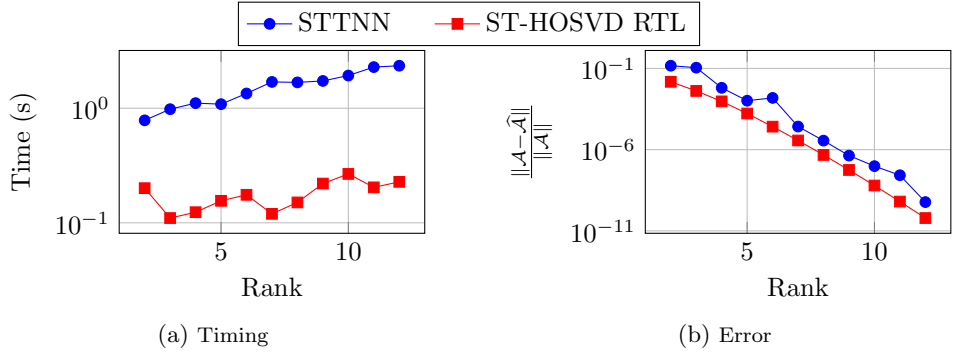

We additionally compare our method to the STTNN approach introduced in~\cite{bucci_randomized_2024}, which also relies on randomized low-rank approximations and employs sequential truncation within each layer. Among the methods evaluated in their work, including Grasedyck's algorithm~\cite[Algorithm~2]{grasedyck_hierarchical_2010} using randomized SVDs, STTNN was reported to be the fastest. To enable a consistent comparison, we implemented our algorithms (\cref{alg:SVD_HTD,alg:RTL_SVD_HTD}) in Matlab as well. The results are shown in \cref{fig:sttnn-strtl}. Across the tested configurations, our method achieves approximately an order-of-magnitude improvement in both runtime and accuracy relative to STTNN.

\subsection{Structural Plate}
\label{sec:plate}
\begin{figure}[tb!]
	\centering
	\begin{subfigure}{0.48\textwidth}
			\centering
			\includegraphics[width=\linewidth]{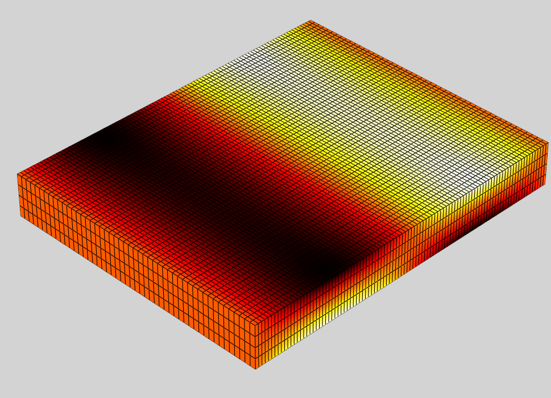}
			\caption{\footnotesize Original X-displacement results: $\tensor{A}[:,:,:,1,1,1,1,1]$}
		\end{subfigure}
	\hfill
		\begin{subfigure}{0.48\textwidth}
			\centering
			\includegraphics[width=\linewidth]{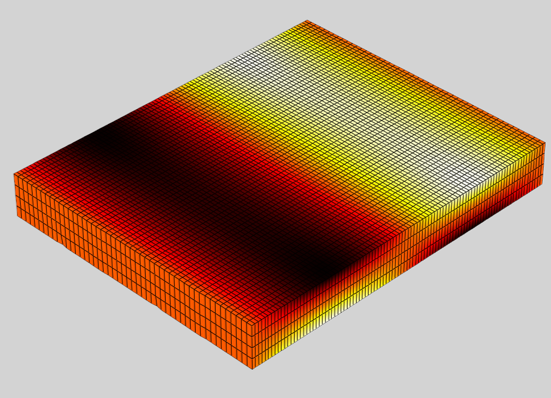}
			\caption{ \footnotesize Projected X-displacement results: $(P_L \circ \dots \circ P_1) (\tensor{A})[:,:,:,1,1,1,1,1]$}
		\end{subfigure}
	\caption{\footnotesize A visualization of the plate finite element mesh, with colors interpolated from the nodal values. The projected version $(P_L \circ \dots \circ P_1) (\tensor{A})$ was computed using the ST-HOVSD RTL algorithm, for this slice $\frac{\|(\tensor{A} - (P_L \circ \dots \circ P_1) (\tensor{A}))[:,:,:,1,1,1,1,1]\|}{\|\tensor{A}[:,:,:,1,1,1,1,1]\|} \approx 0.00765$.}
	\label{fig:plate}
\end{figure}

The first realistic dataset was created by solving the 3D static equilibrium equation $\nabla \cdot \bm{\sigma} + \bm{b} = \bm{0}$ on a plate domain $\Omega=[0,L]\times[-325,325]\times[-T,T]$, where $L\in[800,2000]$~mm is the variable length and $T\in[50,200]$~mm the thickness parameter. The domain was discretised using a structured $100\times50\times5$ mesh of hexahedral finite elements. The problem is solved using the commercial finite element solver Simcenter Nastran (SOL 106), which solves it in its weak form 
\begin{equation} \int_{\Omega} \delta \bm{\varepsilon} : \bm{\sigma} \,\mathrm{d}\Omega = \int_{\Omega} \delta \bm{u} \cdot \bm{b} \,\mathrm{d}\Omega + \int_{\Gamma_t} \delta \bm{u} \cdot \bar{\bm{F}} \,\mathrm{d}\Gamma, \end{equation} where $\bm{\sigma}$ and $\bm{\varepsilon}$ are the stress and strain tensors, respectively, and $:$ denotes the double contraction operator $A : B = \sum_{i,j} A_{i,j} B_{i,j}$. Homogeneous, isotropic linear elasticity $\bm{\sigma}=\mathcal{C}(E):\bm{\varepsilon}$ was assumed, where $\mathcal{C}$ is the fourth-order elasticity tensor, which depends on the young modulus $E \in [70, 550]$~GPa. $\bm{b}$ are the body forces, $\bar{\bm{F}}$ is the applied load, $\delta \bm{u}$ the displacement and $\Gamma_t$ the boundary. The displacement (in X, Y, and Z-direction), the stress (Von Mises), and strain (Von Mises) are the output. The load $\bar{\bm{F}} \in [50, 200]$~kN is equally distributed on a widthwise line in the middle of the plate. All varying parameters are sampled on a 9-point Chebyshev grid, yielding an 8th order tensor of size $100 \times 50 \times 5 \times 5 \times 9 \times 9 \times 9\times 9$ requiring 3.28~GB of storage in single precision. \Cref{fig:plate} shows a simulation result and its reconstruction. 
The dataset tensor is always normalized so that $\|\tensor{A}\| = 1$. 
The 16th-order tensorization of the same data was constructed by splitting its modes into smaller factors, reshaping the tensor to size \((5 \times 4\times 5)\times (5\times 2\times 5)\times (5)\times (3\times 3) \times (3\times 3) \times (3\times 3) \times (3 \times 3)\).
 
\Cref{tab:plate} shows the results for both tensorizations using different methods. While all methods benefit from using ST-HOSVD in general, only LTR shows a clear improvement, almost quartering its execution time. RTL consistently outperforms LTR in speed and allocations, where speedup factors ranging from $10\times$ (against ST-HOSVD) to $39\times$ (against T-HOSVD) 8th order tensor. For RTL, the difference between T-HOSVD and ST-HOSVD is smaller because the binary tree structure allows each Tucker decomposition to only act on 2 modes, limiting the computational benefit of ST-HOSVD. LTR does not benefit much from ST-HOSVD on the 16th order tensorization due to the incompressibility of the leaves. The leaves are incompressible due to the arbitrary tensorization into 16th order, which inhibits compression until the original modes are recombined. The RTL methods also suffer to a lesser degree from the 16th order tensorization due to this incompressibility, although the M-HOSVD method can compress more than RTL leading to improved timing.

Intermediate memory allocations, as measured by \texttt{@benchmark} of BenchmarkTools~\cite{revels_benchmarktoolsjl_2025}, are significantly reduced for all RTL methods compared to LTR, being about $800\times$ for the 16th-order tensor. This reduction stems from both the algorithm and implementation. Algorithmically, the first randomized SVD in the Tucker approximation on the leaf layer of LTR methods will sketch an $N \times N^{d-1}$ matrix down to an $R \times N^{d-1}$ matrix, assuming an $N\times \dots \times N$ input tensor. During RTL methods with a balanced binary tree, the first SVD will sketch an $N^{d/2} \times N^{d/2}$ matrix down to a much smaller $R \times N^{d/2}$. Moreover, the T-HOSVD methods have a higher allocation count than the ST-HOSVD methods due to the implementation's reliance on highly optimized LAPACK calls, which require that all matrices are strided. The matricizations required by T-HOSVD at the leaf layer involve permuting the tensor's original factors, which cannot be realized by simply reinterpreting the linear storage, i.e., $\operatorname{vec}(\tensor{A}_{(k)}) \ne \operatorname{vec}(\tensor{A})$ for $k=2,\dots,d$. Hence, memory must be allocated to store $\tensor{A}_{(k)}$ before performing the sketch. These two reasons contribute to the higher memory allocation of the LTR method with T-HOSVD compared to other methods. 

All methods maintain approximation errors below the prescribed threshold of $10^{-2}$, although the compression ratios can vary. It can be seen, that RTL seems to perform worse than LTR in terms of compression, while M-HOSVD achieves compression closer to LTR. This variance is due to the chosen error distribution and the incompressibility at the leaves. We investigate this in the next subsection.

\begin{table}[tb!]
	\caption{\footnotesize Comparison of T-HOSVD and ST-HOSVD variants on two different tensors. The compression factor is given relative to the uncompressed tensor. The error denotes the relative Frobenius norm of the difference between the original and compressed tensor. The initial sketching size was set to $R^*=200$. The values between parenthesis denote the improvement factors compared to the baseline, being the LTR T-HOSVD variant~\cite[Algorithm 2]{grasedyck_hierarchical_2010}, which is underlined and available in current software. The best results for a tensorization of a dataset are in bold.}
	\label{tab:results}
\begin{subtable}[t]{\linewidth}
	\caption{Plate dataset}
	\label{tab:plate}
	\footnotesize
	\resizebox{0.98\linewidth}{!}{
		\begin{tabular}{@{}lllcccc@{}}
			\toprule
			$\bm{d}$ & \textbf{Tucker} & \textbf{Method} 
			& \textbf{Time [s]} & \textbf{Mem [GiB]} & \textbf{Error} & \textbf{Comp. Factor} \\
			\midrule
			\multirow{6}{*}{8} 
		& \multirow{3}{*}{\underline{T-HOSVD}} 
		& \underline{LTR}   &    $\underline{139.95}\;(1.0\times)$ & $\underline{165.08}\;(1.0\times)$ & \underline{$0.86 \times 10^{-2}$} & $\underline{14192.7}\;(1.0\times)$ \\
		& & RTL       &  $3.66\;(38.2\times)$ & $0.72\;(229.9\times)$ & $\bm{0.76 \times 10^{-2}}$ & $13324.3\;(0.9\times)$ \\
		& & M-HOSVD  & $4.25\;(33.0\times)$ & $0.74\;(223.2\times)$ & $0.77 \times 10^{-2}$ & $14733.5\;(1.0\times)$ \\
		\cmidrule(lr){2-7}
		& \multirow{3}{*}{ST-HOSVD} 
		& LTR       &    $37.79\;(3.7\times)$ & $42.70\;(3.9\times)$ & $0.86 \times 10^{-2}$ & $14192.7\;(1.0\times)$ \\
		& & RTL       &  $\bm{3.57\;(39.2\times)}$ & $\bm{0.65\;(254.7\times)}$ & $0.77 \times 10^{-2}$ & $13630.8\;(1.0\times)$ \\
		& & M-HOSVD  & $3.62\;(38.7\times)$ & $0.67\;(246.8\times)$ & $0.78 \times 10^{-2}$ & \bm{$15080.8\;(1.1\times)$} \\
			\midrule 
			\multirow{6}{*}{16} 
			& \multirow{3}{*}{\underline{T-HOSVD}} 
			& \underline{LTR}   &    $\underline{365.51}\;(1.0\times)$ & $\underline{572.44}\;(1.0\times)$ & \underline{$0.94 \times 10^{-2}$} & $\underline{31046.5}\;(1.0\times)$ \\
			& & RTL       &  $8.73\;(41.9\times)$ & $0.70\;(817.2\times)$ & $\bm{0.47 \times 10^{-2}}$ & $15772.5\;(0.5\times)$ \\
			& & M-HOSVD  & $\bm{8.50\;(43.0\times)}$ & $0.73\;(786.6\times)$ & $0.79 \times 10^{-2}$ & $25137.2\;(0.8\times)$ \\
			\cmidrule(lr){2-7}
			& \multirow{3}{*}{ST-HOSVD} 
			& LTR       &    $258.95\;(1.4\times)$ & $449.90\;(1.3\times)$ & $0.94 \times 10^{-2}$ & $\bm{31630.9\;(1.0\times)}$ \\
			& & RTL       &  $8.60\;(42.5\times)$ & $\bm{0.64\;(892.7\times)}$ & $0.53 \times 10^{-2}$ & $16904.9\;(0.5\times)$ \\
			& & M-HOSVD  & $8.80\;(41.5\times)$ & $0.67\;(856.3\times)$ & $0.80 \times 10^{-2}$ & $24120.6\;(0.8\times)$ \\
			\bottomrule
		\end{tabular}
	}	
\end{subtable}
\medskip
	\begin{subtable}[b]{\linewidth}
	\caption{``Small'' Miranda dataset}
	\label{tab:miranda}
	\footnotesize
	\resizebox{0.98\linewidth}{!}{
	\begin{tabular}{@{}lllcccc@{}}
			\toprule
			$\bm{d}$ & \textbf{Tucker} & \textbf{Method} 
			& \textbf{Time [s]} & \textbf{Mem [GiB]} & \textbf{Error} & \textbf{Comp. Factor} \\
			\midrule
			\multirow{6}{*}{4}
			& \multirow{3}{*}{\underline{T-HOSVD}} 
			& \underline{LTR}   &    $\underline{20.47}\;(1.0\times)$ & $\underline{22.03}\;(1.0\times)$ & \underline{$0.99 \times 10^{-2}$} & $\bm{\underline{101.3}\;(1.0\times)}$ \\
			& & RTL       &  $2.58\;(7.9\times)$ & $1.78\;(12.4\times)$ & $0.94 \times 10^{-2}$ & $93.2\;(0.9\times)$ \\
			& & M-HOSVD  & $2.80\;(7.3\times)$ & $2.19\;(10.0\times)$ & $\bm{0.80 \times 10^{-2}}$ & $92.8\;(0.9\times)$ \\
			\cmidrule(lr){2-7}
			& \multirow{3}{*}{ST-HOSVD} 
			& LTR       &    $6.73\;(3.0\times)$ & $6.16\;(3.6\times)$ & $0.99 \times 10^{-2}$ & $\bm{101.3\;(1.0\times)}$ \\
			& & RTL       &  $\bm{2.21\;(9.3\times)}$ & $\bm{1.32\;(16.7\times)}$ & $0.95 \times 10^{-2}$ & $93.8\;(0.9\times)$ \\
			& & M-HOSVD  & $2.40\;(8.5\times)$ & $1.71\;(12.9\times)$ & $0.81 \times 10^{-2}$ & $92.8\;(0.9\times)$ \\
			\midrule 
			\multirow{6}{*}{16} 
			& \multirow{3}{*}{\underline{T-HOSVD}} 
			& \underline{LTR}   &    $\underline{103.99}\;(1.0\times)$ & $\underline{192.18}\;(1.0\times)$ & \underline{$1.00 \times 10^{-2}$} & $\bm{\underline{99.4}\;(1.0\times)}$ \\
			& & RTL       &  $5.97\;(17.4\times)$ & $2.90\;(66.3\times)$ & $\bm{0.59 \times 10^{-2}}$ & $60.7\;(0.6\times)$ \\
			& & M-HOSVD  & $7.18\;(14.5\times)$ & $4.34\;(44.3\times)$ & $0.79 \times 10^{-2}$ & $78.3\;(0.8\times)$ \\
			\cmidrule(lr){2-7}
			& \multirow{3}{*}{ST-HOSVD} 
			& LTR       &    $71.03\;(1.5\times)$ & $148.26\;(1.3\times)$ & $1.00 \times 10^{-2}$ & $90.9\;(0.9\times)$ \\
			& & RTL       &  $\bm{5.58\;(18.6\times)}$ & $\bm{2.38\;(80.8\times)}$ & $0.62 \times 10^{-2}$ & $61.1\;(0.6\times)$ \\
			& & M-HOSVD  & $6.48\;(16.1\times)$ & $3.79\;(50.7\times)$ & $0.80 \times 10^{-2}$ & $77.6\;(0.8\times)$ \\
			\bottomrule
	\end{tabular}
}
\end{subtable}
\end{table}

\subsection{Miranda CFD dataset}
\label{sec:miranda}
\begin{figure}[tb]
	\centering
	\begin{subfigure}{0.48\textwidth}
		\centering
		\includegraphics[width=\linewidth]{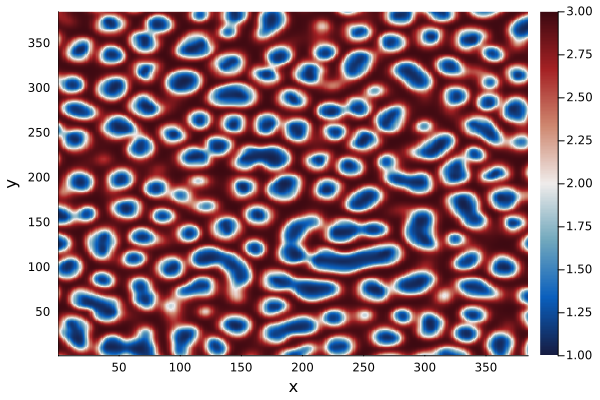}
		\caption{\footnotesize 145th z-slice of the original tensor: $\tensor{A}[:,:,145,1]$}
		\label{fig:miranda_origin}
	\end{subfigure}
	\hfill
	\begin{subfigure}{0.48\textwidth}
		\centering
		\includegraphics[width=\linewidth]{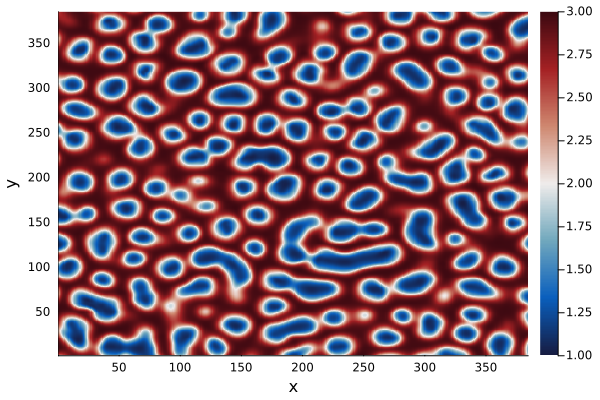}
		\caption{\footnotesize 145th z-slice of the projected tensor: $(P_L \circ \dots \circ P_1) (\tensor{A})[:,:,145,1]$}
		\label{fig:miranda_recon}
	\end{subfigure}
	\caption{\footnotesize The density on the $xy$-domain of the ``Small'' Miranda dataset from \cite{zhao_sdrbench_2020} at $z$-index 145. The projected version $(P_L \circ \dots \circ P_1) (\tensor{A})$ was computed using the ST-HOVSD RTL algorithm.}
	\label{fig:miranda}
\end{figure}

The last dataset is the ``small'' Miranda computational fluid dynamics dataset~\cite{cook_tera-scalable_2005} from the Scientific Data Reduction Benchmark (SDRBench\footnote{https://sdrbench.github.io}) of Zhao et al~\cite{zhao_sdrbench_2020}. According to the metadata, the data was subsampled from a simulation of a Rayleigh--Taylor instability on a $3,072^3$ uniform grid, resulting in a $384 \times 384 \times 256 \times 7$ tensor. The first three factors correspond to the spatial domain and the last mode contains the 7 output variables:\footnote{The author of the dataset states the pressure, viscosity and diffusivity are mislabeled, but we keep them in accordance with the description of the dataset metadata.} density, pressure, viscosity, diffusivity, and 3 velocity components~\cite{cook_tera-scalable_2005}. The fibers corresponding to each of these outputs was normalized individually before normalizing the full tensor.
\Cref{fig:miranda_origin} shows a slice of the original tensor in the $z$-dimension at index 145 for the density, while \cref{fig:miranda_recon} is the same slice but taken from the approximated tensor as computed by the ST-HOSVD RTL method.
The experimental results of the decomposition of the Miranda tensor can be seen in \cref{tab:miranda}.
Most of the remarks from \cref{sec:plate} hold in this case as well. The difference between ST-HOSVD and T-HOSVD is more pronounced, leading to speedup factors of $1.5\times$ and $3\times$ for the LTR case depending on the tensorization. The 16-order tensorization benefits less from ST-HOSVD because the leaves are incompressible due to the arbitrary tensorization, negating the benefit of using the ST-HOSVD.
RTL methods again outperform the LTR variants in execution time and memory allocation counts. The speedup factor for the original dimension varies between $3\times$ and $9\times$ depending on whether ST-HOSVD or T-HOSVD LTR is used, while the 16th order tensorization achieves even greater speedups.

When comparing the compression factors, LTR methods outperform the RTL versions. Especially the simple RTL methods, without the SVD-HTD correction, seem to suffer from the order-16 tensorization. We attribute this to the way the error estimation is used to guide truncation and the observed incompressibility of the leaves. Due to the redistribution of the leftover error ``budget,'' nodes at the end of the algorithm are allowed to make a larger error than the initial ones. This explains why LTR achieves a better compression factor, as it can use the unspent error budget of the leaves, at the nodes closer to the root, while RTL is allocating a large allowed error for the leaves, which cannot be used effectively. This effect is evident in the smaller error of RTL compared to LTR. On top of this, the simple RTL method overestimates the error it makes since it does not compute the exact singular values. SVD\_RTL suffers less from this.
 A possible alleviation of this problem, in a practical setting, would be to apply weights to the allowed error at each node, so that nodes higher up in the tree are allowed to make a larger error than the nodes near the leaves. As the possible rank and number of singular values near the root is larger, these nodes have the possibility of larger compression rates. One could also run the truncation algorithm again on the HTD representation, which is comparatively cheap~\cite{hackbusch_tensor_2019}.

\section{Conclusions}
\label{sec:conclusions}
This paper demonstrates that when randomized SVDs are used, the LTR direction used in the current state-of-the-art HTD libraries can be one order of magnitude less efficient than the RTL direction. While the quasi-optimal error bound of RTL is larger than that of LTR, due to the non-orthogonality of intra-layer truncation errors, this difference was not pronounced in practical cases. Experimental results also support the theoretical computational complexity results, showing speedup factors of RTL over LTR in the range of $3\times$ to $40\times$ and total memory allocation reductions of at least $10 \times$, while obtaining compression factors similar to the most commonly used HTD method of Grasedyck~\cite[Algorithm 2]{grasedyck_hierarchical_2010}. 

The computational complexity and error bound of the RTL algorithm show a clear preference for square matrix flattenings. It would thus be preferable to compute the decomposition using a balanced tree and afterwards convert this tree structure into a desired form, although this is the subject of future research. Better error distribution across layers to achieve a given error tolerance could also be investigated.

\section*{Acknowledgments}
KU Leuven's licensed Microsoft Copilot was exclusively used to create \LaTeX{} code templates for tables and figures, and to check for grammatical and spelling errors. The Siemens' licensed GitHub Copilot created the 3D visualization code for the plate use case. All other content was conceived and written by humans and the authors assume responsibility for all content.  
We also acknowledge our used dataset coming from the SDRBench Miranda Dataset, the DOE NNSA ECP project, and the ECP CODAR project. We thank Alberto Bucci and Sebastian Kr{\"a}mer for their insightful comments at METT XI and for pointing out that the root-independent normal form had already been described in~\cite{Grasedyck2024}.

\bibliographystyle{siamplain}
\bibliography{references.bib}

\appendix
\section{Proof of \cref{alg:SVD_HTD}}
\label{app:svd_htd_proof}
\begin{proposition}
		\Cref{alg:SVD_HTD} transforms a general HTD of $\tensor{A} \in \mathrm{V}_1 \otimes \dots \otimes \mathrm{V}_d$ with dimension tree $\mathcal{T}$ into SVD-HTD form of $\tensor{A}$.
		\begin{proof}
			After performing an LTR QR orthonormalization in line 1 using~\cite[Algorithm 3]{grasedyck_hierarchical_2010}, every non-root node $t$ stores a $B_t$ with orthonormal columns. Then, by~\cite{hackbusch_tensor_2019}, all frames $F_t$ defined recursively by \cref{eq:rec_def}, with $F_k = B_k$ for all leaves $k \in \mathcal{T}_L$, also have orthonormal columns, except at the root. For the root node, we have $F_1 = \operatorname{vec}(\tensor{A})$ and therefore $\|F_1\| = \|\tensor{A}\|$, so the SVD-HTD property is already satisfied at the root layer. We show that one iteration sweep preserves the SVD-HTD property for the parents while establishing it for its children
			
			Consider one iteration step. In lines~5--6 a non-truncating HOSVD is applied to $\tensor{B}_t$ and its children's transfer tensors $B_{t_i}$ are transformed accordingly. The update $B_{t_j} \gets B_{t_j}U_{t_j} \Sigma_{t_j}$ induces the frame update $F_{t_j} \gets F_{t_j}U_{t_j} \Sigma_{t_j}$ and therefore
			\[
			\begin{aligned}
				F_t &\gets (F_{t_1}U_{t_1} \Sigma_{t_1} \otimes \dots \otimes F_{t_m}U_{t_m} \Sigma_{t_m})(\Sigma_{t_1}^{-1} U_{t_1}^\top \otimes \dots \otimes \Sigma_{t_m}^{-1} U_{t_m}^\top)  \tensor{B}_t \\
				&= (F_{t_1}\otimes \dots \otimes F_{t_m})  \tensor{B}_t = F_t,
			\end{aligned}
			\]
			where we assume compact SVDs are used so that $\Sigma^{-1}$ is properly defined.
			Thus, $F_t$ is unchanged and the SVD-HTD property is preserved at the parent. Since the SVD-HTD property was satisfied at the parent $t$, \cref{cor:HOSVD_child} ensures the $\Sigma_{t_i}$ are the singular values of $\tensor{A}_{t_i}$. Additionally, since each $F_{t_i}$ contains orthonormal columns, it is easily verified that $F_{t_i}U_{t_i}$ are the corresponding left singular vectors. Thus the children are put into SVD-HTD form by line 6, concluding the proof.
		\end{proof}
	\end{proposition}
\section{Discretized Helmholtz equation}
\label{sec:helm}
\begin{figure}[tb!]
	\centering
	\begin{subfigure}{0.48\textwidth}
			\centering
			\includegraphics[width=\linewidth]{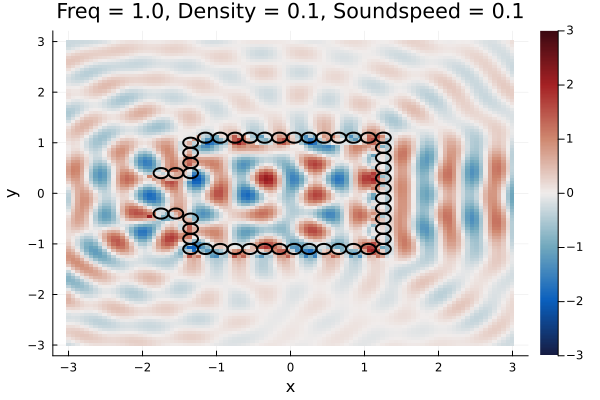}
			\caption{\footnotesize First frontal slice of the original tensor: $\tensor{A}[:,:,1,1,1]$}
		\end{subfigure}
	\hfill
		\begin{subfigure}{0.48\textwidth}
			\centering
			\includegraphics[width=\linewidth]{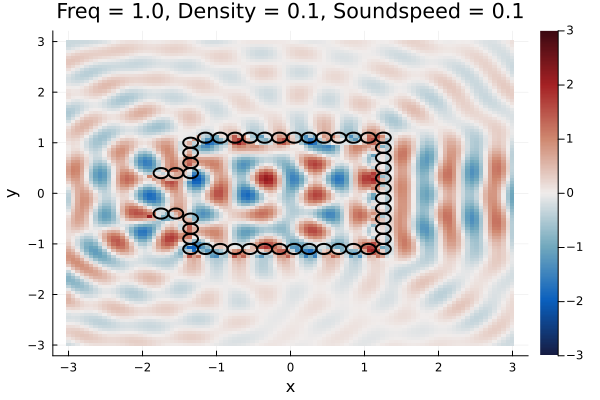}
			\caption{ \footnotesize First frontal slice of the projected tensor: $(P_L \circ \dots \circ P_1) (\tensor{A})[:,:,1,1,1]$}
		\end{subfigure}
	\caption{\footnotesize The real part of the solution to the Helmholtz resonator as computed by \cite{gower_multiplescatteringjl_2025}. The projected version $(P_L \circ \dots \circ P_1) (\tensor{A})$ was computed using the ST-HOVSD RTL algorithm.}
	\label{fig:helm}
\end{figure}

Another realistic dataset originates from solving the 2D Helmholtz equation $\nabla^2 u(x,y,\omega) + \frac{\omega}{c}^2u(x,y,\omega) = 0$ on a uniform \(100\times 100\) spatial $xy$-grid in $ [-3,3] $ $\times [-3,3]$ (in m) for 50 angular frequencies $\omega$ in $[1,5]$ (in Hz), 20 medium densities $\rho$ in $[0.1,2]$ (in kg/m$^3$), and 20 sound speeds $c$ in $[0.1,2]$ (in m/s). The Helmholtz resonator example from the MultipleScattering.jl library \cite{gower_multiplescatteringjl_2025} was used for this test case, where the resonator elements have a density of $10$ kg/m$^3$ and the speed of sound inside the element was set to $5$ m/s. The difference in density $\rho$ between the surrounding medium and the resonator influences the reflection of sound. A simulation result and its reconstruction, computed by contracting the full dimension tree, is shown in \cref{fig:helm}. 

Storing all parameter variations yields an order-5 tensor of size \(100 \times 100 \times 50 \times 20 \times 20\) requiring 1.6 GB of storage. The dataset tensor is always normalized so that $\|\tensor{A}\| = 1$. 
The 16th-order tensorization of the same data was constructed by splitting its modes into smaller factors, producing a tensor of size \((5 \times 2\times 5\times 2)\times (5\times 2\times 5\times 2)\times (5\times 2\times 5)\times (2\times 2\times 5)\times (4\times 5)\) with the same number of elements. 

\begin{table}[tb!]
	\caption{\footnotesize Comparison of T-HOSVD and ST-HOSVD variants on two different tensors. The compression factor is given relative to the uncompressed tensor. The error denotes the relative Frobenius norm of the difference between the original and compressed tensor. The initial sketching size was set to $R^*=200$. The values between parenthesis denote the improvement factors compared to the baseline, being the LTR T-HOSVD variant~\cite[Algorithm 2]{grasedyck_hierarchical_2010}, which is underlined and available in current software. The best results for a tensorization of a dataset are in bold.}
	\label{tab:extra_results}
		\begin{subtable}[t]{\linewidth}
				\caption{Helmholtz dataset}
				\label{tab:helmholtz}
				\footnotesize
				\resizebox{0.98\linewidth}{!}{
					\begin{tabular}{@{}lllcccc@{}}
								\toprule
								$\bm{d}$ & \textbf{Tucker} & \textbf{Method} 
								& \textbf{Time [s]} & \textbf{Mem [GiB]} & \textbf{Error} & \textbf{Comp. Factor} \\
								\midrule
								\multirow{6}{*}{5} 
				& \multirow{3}{*}{\underline{T-HOSVD}} 
				& \underline{LTR}   &    $\underline{31.55}\;(1.0\times)$ & $\underline{41.82}\;(1.0\times)$ & \underline{$0.98 \times 10^{-2}$} & $\bm{\underline{11.2}\;(1.0\times)}$ \\
				& & RTL       &  $11.49\;(2.7\times)$ & $11.29\;(3.7\times)$ & $\bm{0.84 \times 10^{-2}}$ & $9.9\;(0.9\times)$ \\
				& & M-HOSVD  & $14.93\;(2.1\times)$ & $15.44\;(2.7\times)$ & $0.86 \times 10^{-2}$ & $10.2\;(0.9\times)$ \\
				\cmidrule(lr){2-7}
				& \multirow{3}{*}{ST-HOSVD} 
				& LTR       &    $24.53\;(1.3\times)$ & $30.69\;(1.4\times)$ & $0.99 \times 10^{-2}$ & $\bm{11.2\;(1.0\times)}$ \\
				& & RTL       &  $\bm{9.24\;(3.4\times)}$ & $\bm{8.12\;(5.1\times)}$ & $0.92 \times 10^{-2}$ & $9.9\;(0.9\times)$ \\
				& & M-HOSVD  & $12.90\;(2.4\times)$ & $12.11\;(3.5\times)$ & $0.91 \times 10^{-2}$ & $10.1\;(0.9\times)$ \\
								\midrule 
								\multirow{6}{*}{16} 
				& \multirow{3}{*}{\underline{T-HOSVD}} 
				& \underline{LTR}   &    $\underline{101.83}\;(1.0\times)$ & $\underline{176.29}\;(1.0\times)$ & \underline{$1.00 \times 10^{-2}$} & $\bm{\underline{9.1}\;(1.0\times)}$ \\
				& & RTL       &  $6.97\;(14.6\times)$ & $4.52\;(39.0\times)$ & $\bm{0.36 \times 10^{-2}}$ & $6.7\;(0.7\times)$ \\
				& & M-HOSVD  & $15.81\;(6.4\times)$ & $13.67\;(12.9\times)$ & $0.80 \times 10^{-2}$ & $8.5\;(0.9\times)$ \\
				\cmidrule(lr){2-7}
				& \multirow{3}{*}{ST-HOSVD} 
				& LTR       &    $95.38\;(1.1\times)$ & $167.14\;(1.1\times)$ & $1.00 \times 10^{-2}$ & $8.9\;(1.0\times)$ \\
				& & RTL       &  $\bm{6.80\;(15.0\times)}$ & $\bm{4.20\;(41.9\times)}$ & $0.40 \times 10^{-2}$ & $6.8\;(0.7\times)$ \\
				& & M-HOSVD  & $15.62\;(6.5\times)$ & $13.19\;(13.4\times)$ & $0.80 \times 10^{-2}$ & $8.4\;(0.9\times)$ \\
								\bottomrule
						\end{tabular}
				}	
			\end{subtable}
\end{table}

\Cref{tab:helmholtz} shows the results for both tensorizations using different methods.
As seen for the other datasets, the ST-HOSVD variants consistently achieve lower computational times compared to their T-HOSVD counterparts, although the difference is rather limited.

Interestingly, for the RTL methods, the order-16 tensorization yields a substantial reduction ($\approx 1.5 \times$) in execution time relative to the order-5 case. The increased order leads to more uniform mode sizes, which results in more square-like matricizations during the algorithm, when a balanced dimension tree is used. As indicated in \cref{sec:motivation,sec:complexity}, square matricizations can reduce the computational complexity. 

RTL still consistently outperforms LTR, where speedup factors ranging from $3 \times$ (against ST-HOSVD) to $15\times$ (against HOSVD) can be seen. For RTL, the difference between T-HOSVD and ST-HOSVD is smaller, because the binary tree structure allows each Tucker decomposition to only act on 2 modes, limiting the computational benefit of ST-HOSVD. The reason why LTR does not benefit much from ST-HOSVD is due to the near incompressibility of the leaves. Compression is realized only in the layers closer to the root in these experiments.

The difference between the RTL variants can be quite large. This is due to the rather large hierarchical ranks that can make the nodewise SVD-HTD transformation costly. For example, the HTD of the order-16 tensorization had hierarchical ranks of about $1100$ in the children of the root. Despite the implications of this higher rank in light of the complexity analysis in \cref{sec:complexity}, RTL is still much faster than LTR. 

\end{document}